\documentclass[reqno]{amsart}
\usepackage{amsfonts}
\usepackage{amssymb}
\usepackage[T1]{fontenc}
\usepackage{ae}
\usepackage{hyperref}

\title{Fixed Point Sets For Permutation Modules}
\author{Peter Collings}
\begin{document}
\maketitle
  
\newcounter{Subsection}[section] \newcounter{chapterno}
\setcounter{chapterno}{1} \newcounter{lastplace}
\newcommand{\place}{{\bf \arabic{chapterno}.\arabic{lastplace}.}}
\newcommand{\Subsection}[2]{\addtocounter{Subsection}{1}
  \setcounter{lastplace}{} \newpage
  #2
    \begin{center} \LARGE{\arabic{section}.\arabic{Subsection}}\ \ \ \ \ {\LARGE\texttt {#1}} \end{center}   
    \addtocounter{lastplace}{1}} \newenvironment{theorm}[2] {
    \noindent \medskip \\  {\bf { Theorem}} {\large{\place}} {\it { #1}}
    \noindent \medskip \\ {\bf Proof.} #2$\! \! \! \! \! \qed$}{\noindent 
    \addtocounter{lastplace}{1}} \newenvironment{lemma}[2] {
    \noindent \medskip \\  {\bf { Lemma}} {\large{\place}} {\it { #1}}
    \noindent \medskip \\ {\bf Proof.} #2$\! \! \! \! \! \qed$}{\noindent 
    \addtocounter{lastplace}{1}} \newenvironment{corollary}[2] {
    \noindent \medskip \\  {\bf { Corollary}} {\large{\place}} {\it { #1}}
    \noindent \medskip \\ {\bf Proof.} #2$\! \! \! \! \! \qed$}{\noindent 
    \addtocounter{lastplace}{1}} \newenvironment{remark}[1] {
    \noindent \medskip \\  {\bf { Remark}} {\large{\place}} {#1}}
  {\noindent \addtocounter{lastplace}{1}} \newenvironment{block}[2] {
    \noindent \medskip \\  
    #2{\large{\place}} { #1}} {\noindent \addtocounter{lastplace}{1}}
  
  \newenvironment{nblock}[2] { \setcounter{lastplace}{1}
     \noindent \medskip \\ \begin{center}
       {\large{\arabic{chapterno}.}} {\MakeUppercase {#2} } \newline
     \end{center} { #1} }

   \newenvironment{fact}[1] {
    \noindent \medskip \\  
    {\bf { Fact}} {\large{\place}}{\it { #1}}} {\noindent
    \addtocounter{lastplace}{1}} \newenvironment{definition}[1] {
    \noindent \medskip \\  
    {\bf { Definition}} {\large{\place}} #1} {\noindent
    \addtocounter{lastplace}{1}} \newenvironment{convention}[1] {
    \noindent \medskip \\  
    {\bf { Convention}} {\large{\place}} #1} {\noindent
    \addtocounter{lastplace}{1}} \newenvironment{example}[1] {
    \noindent \medskip \\  {\bf { Example}} {\large{\place}} {#1}}
  {\noindent \addtocounter{lastplace}{1}}
  
  \newcommand{\zerostart}{\addtocounter{Subsection}{-1}}
  \newcommand{\zeroblock}{\addtocounter{lastplace}{-1}}
  \setcounter{lastplace}{1} \newcommand{\mylabel}[1] { \newcounter{#1}
    \setcounter{#1}{\arabic{lastplace}} \newcounter{#1_cn}
    \setcounter{#1_cn}{\arabic{chapterno} } } \newcommand{\partn}{
    {\tiny{ | \! - } } }

  \newcommand{\arabicc}[2]{(\arabic{#1_cn}.\arabic{#1})}
  
  \newcommand{\Fix}{{\rm Fix}} \newcommand{\Stab}{{\rm Stab}}
  \newcommand{\supp}{{\rm supp}}
  \newcommand{\Sym}{{\rm Sym}}
  \newcommand{\Syl}{{\rm Syl}}

\numberwithin{equation}{subsection}
\baselineskip=12pt
\textwidth=15cm
\textheight=21cm

  \vspace{0cm}
\markboth{Preliminaries}{Preliminaries}
Let $ A = k G $, the group algebra of some finite group where the
characteristic of the field $ k $ divides $ | G | $. In contrast to
working over the complex field, the $ k G $-modules are not usually
semisimple. If a Sylow $ p $-subgroup of $ G $ is not cyclic then
there are infinitely many indecomposable $ k G $-modules, and we
usually enjoy little control over the category of such modules. It is
therefore an important problem to find classes of modules which may be
expressed as a sum of a not very great number of indecomposables, and
to understand the structure of these indecomposables.
        \par
        \vspace{5mm}    Permutation $ k G$-modules, and their
        indecomposable summands (called \linebreak $ p $-permutation
        modules), provide one such class. It is only necessary  to
        consider transitive permutation actions and for a fixed group
        this class is a finite set. Furthermore each indecomposable
        summand $ M $ determines three invariants, namely a vertex, a
        quotient of a subgroup of $ G $ and a projective
        indecomposable  $ p $-permutation module for this quotient;
        and these three invariants characterise $ M $ up to
        isomorphism. These invariants   were introduced by Green
        \cite{gre58} and Brauer         and adapted by Puig and, in
        the case of permutation modules, by Brou{\'e} \cite{bro85}.
        \par 
        Given an indecomposable $ k G $-module $ M $, a vertex of $ M
        $ is a subgroup $ H \leq G $ minimal subject to satisfying the
        condition that there be a $ k H $-module $ L $ for which $ M $
        is a component of $ {\rm Ind} ^ G _ H ( L ) $. This vertex is
        unique  up to conjugation in $ G $. Any Sylow $ p $-subgroup
        of $ G $ satisfies this condition and so vertices are always $
        p $-subgroups.  The $ k H $-module      $ L $ is called a
        source of $ M $ and is unique up to conjugation in $ N_G ( H )
        $. However, $ M $ is a $ p $-permutation        module if and
        only if $ M $ is a component of a permutation module and this
        happens exactly when $ L $ is the trivial       $ k H
        $-module. Hence $ p $-permutation modules are also known as
        trivial source modules ({\it cf.} \cite[p216]{the95}). \par
        \par 
        The Brauer correspondent $ M ^ H $ of $ M $ is an
        indecomposable $ k ( N_G ( H ) / H ) $-module. When $ M $ is a
        component of a permutation module $ k \Omega $ then $ M ^ H $
        is a component of $ k {\rm Fix}_{\Omega} H $, hence is a $ p
        $-permutation module. Brou\'e's correspondence (1.1) asserts a
        $ 1- 1 $ correspondence between the projective components of $
        k {\rm Fix}_{\Omega} H $ and the components of $ M $ with
        vertex $ H $.
        \par \vspace{5mm}
        In Theorem (2.11) we will seek to phrase Brou\'e's
        correspondence  in terms not of vertices but of {\it fixed
          point sets} (defined in (2.4)), a notion which we will see
        is equivalent to that of vertices.
        \par 
        In Section 3 and onwards we shall let $ n $ be a positive
        integer and utilize fixed point sets in the situation of
        permutation modules afforded by the action by conjugation of
        the symmetric group $ {\rm Sym} ( q n ) $       on its
        conjugacy class of fixed point free elements of order $ q $.
        Our aim will be to determine as far as possible the fixed
        point sets of the components involved. At the same time
        we will show how the corresponding
        vertices may be computed eg Lemma (4.5) and  Lemma (4.6), 
        and (5.4).
        \par
        In Theorem (7.18) we determine the general form of a fixed
        point set. From this theorem a list can be derived of
        possible fixed point sets in terms not just of irreducible
        sets but of the even more basic transitive sets. 
        For a given $ p$ and $ q $ we need the further 
        information of: $ (1) $
        those values of $ n $ for which our permutation module admits
        a projective    summand and $ (2) $ the numbers $ \kappa ( X
        ) $ (defined in (7.15)) for irreducible fixed point
        sets $ X $. These missing values are not investigated
        in the present paper.
        
        \par When $ p = q = 2 $ the problem becomes to
        determine the fixed point sets, or the vertices
        of the components of, the
        permutation module afforded by the action of $ \Sym ( 2 n ) $
        on its conjugacy class $ \Xi_{ 2 n } ^ 2 $ of fixed point free
        involutions. This permutation action has been studied before;
        as an example the ordinary character it affords is known, and
        we refer to \cite{ing90} for a demonstration of this and
        associated results. In the situation $ p =q = 2 $ the values
        labelled $ (1) $ and $ (2) $ above can be obtained, leading to
        a precise list of vertices. In the author's PhD thesis
        \cite{phd04} this precise list is       constructed, along
        with an explicit determination of the number of components and
        the Brauer quotients. These results     have been collected by
        the author elsewhere, along with a determination of the Brauer
        characters of the components, again     using a result on
        ordinary characters contained in \cite{ing90}.
        \par The tool        
        to achieve the above is to view the
        natural embeddings of the cartesian product of any two
        symmetric groups into a third as a type of binary multiplication
        which will allow us to construct new fixed point sets
        from existing ones.
        The
        binary multiplication even allows for
        an essentially unique decomposition into irreducibles,
        with products of coprime
        fixed point sets being easy to control.
        Furthermore a unary multiplication is
        provided by the diagonal subset of the  cartesian product.      
        
        \markright{Preliminaries.}
\begin {nblock}
  {      Let $ G $ be a finite group and let $ P $ be a Sylow $ p
    $-subgroup of $ G $. A $ k G $-module has a $ P $-invariant basis
    if and only if it is a component of some permutation $ k
    G$-module, and such a module is called a $ p $-permutation module.
    The main tool for identifying the vertices of the components of
    permutation modules is the following correspondence due to Brou\'e
    [1, (3.2)]. The theorem following it (see [6, (27.7)]) is another
    useful characterisation of these vertices.
}{Preliminaries.}%
\end{nblock}
\mylabel{BCorr}
\begin{block}
  {      {\it Let $ \Omega $ be a permutation $ G $-space. There is a
      multiplicity-preserving bijective correspondence between the
      components of the permutation $ k G$-module $ k \Omega $ with
      vertex $ Q \leq G $ and the projective components of the
      permutation $ k ( N_G ( Q ) / Q) $-module $ k \Fix_{ \Omega } (
      Q ) $.}    } {{\bf Brou\' e Correspondence} }
\end {block}

\begin{block}
  {     {\it Let $ P $ be a Sylow $ p $-subgroup of $ G $.  An
      indecomposable $ k G $-module $ M $ with a $ P $-invariant basis
      $ X $ has vertex conjugate to $ Q \leq P $ if and only if: $ Q $
      is a maximal subgroup of $ P $ subject to fixing an element of $
      X $.}
    
  } {{\bf Theorem} }
\end{block}
\par Recall that any transitive permutation
$ k P $-module is indecomposable (eg \cite[(27.1)]{the95}).  Thus the
regular module $ k P $ is indecomposable and the only indecomposable
projective $ k P $-module is $ k P $.
\mylabel{lemma:RegularProjective}
\begin {lemma}
  {     Let $ M $ be an indecomposable permutation $ P $-module with a $ P
    $-invariant basis $ X $ and let $ x \in X $.  Then $ \Stab_P ( x )
    $ is a vertex of $ M $. In particular $ P $ acts regularly on $ X
    $ if and only if $ M $ is projective.  }  {  First of all we note
    that $ P $ must act transitively on $ X $ because $ M $ is
    indecomposable.  A vertex of $ M $ is any subgroup $ Q $ of $ P $
    maximal subject to fixing an element of $ X $. By taking a
    conjugate of $ Q $ in $ P $ we may assume that $ Q $ fixes the
    element $ x \in X $. Thus $ Q \leq \Stab_P ( x ) $ and the
    maximality of $ Q $ gives $ Q = \Stab_P ( x ) $. Now $ M $ is
    projective if and only if $ Q $ is trivial, and since $ Q =
    \Stab_P ( x ) $ this is the case if and only if $ P $ acts
    regularly on $ X $.  }
\end {lemma}
\vspace*{1mm}
\newline The following is a standard result.
\mylabel{lemma:ProjectiveModules}
\begin {lemma}
  {     {\bf Projective Modules.} Let $ M $ be a $ k G $-module and
    $ P $ a Sylow $ p $-subgroup of $ G $.  Then
    $ M $ is projective if and only if $ {\rm Res} ^ G _ P ( M ) $ is
    projective.
}
{ 
}
\end {lemma}
\begin {lemma}
  {     A $ p $-permutation $ k G $-module $ M $ with $ P $-invariant
    basis $ X $ is projective if and only if each $ P $-orbit on $ X $
    is regular.  }  {   Let $ X_1, X_2 , \dots , X_s $ be the $ P
    $-orbits on $ X $. Then $ {\rm Res} ^ G _ P ( M ) = k X_1 \oplus k
    X_2 \oplus \dots \oplus k X_s $ is a decomposition of $ {\rm Res}
    ^ G _ P ( M ) $ into indecomposable summands.  Now by Lemma
    \arabicc{lemma:ProjectiveModules}{} $ M $ is projective if and
    only if $ {\rm Res} ^ G _ P ( M ) $ is projective, and by
    \linebreak Lemma \arabicc{lemma:RegularProjective}{} this happens
    if and only if $ P $ acts regularly on each $ X_i $.  }
\end {lemma}
\vspace*{0.1cm} \par The following is well known.
\begin {block}
  {     {\it Let $ G $ and $ H $ be finite groups. Let $ \{ S_i
      \}_{i=1}^s $ be a complete set of pairwise non-isomorphic
      indecomposable projective $ k G $-modules and $ \{ T_i
      \}_{i=1}^t $ a complete set of pairwise non-isomorphic
      indecomposable projective $ k H $-modules. Then $ \{ S_i \otimes
      T_j \}_{ 1 \leq i \leq s , 1 \leq j \leq t } $ is a complete set
      of pairwise non-isomorphic indecomposable projective $ k ( G
      \times H ) $-modules.}     } {{\bf Proposition} }
\end {block}
\mylabel{corollary:ProdHasProj}
\begin {corollary}
  {     If $ G $ is a finite group and $ M $ is a $ k G $-module let $
    np \ (M) $ be the number of indecomposable projective summands of
    $ M $ counting multiplicities when $ M$ is expressed as a direct
    sum of indecomposabl;e modules.  Then if $ G_1 $ and $ G_2 $ are
    finite groups and $ M_i $ is a $ k G_i $-module for $ i = 1 ,2 $
    we have
    $$
    np \ ( M_1 \otimes M_2) = np \ (M_1) \cdot np \ (M_2) .$$
    \vspace*{-1cm} } { }
\end {corollary}
\markright{Wreath Products}
\mylabel{section:WreathProd}
\begin {block}
  {      Let $ M $ be a $ k G $-module.  Let $ a $ be a positive
    integer and let $ G ^ a := G \times G \times \dots \times G $ be
    the direct product of $ a $ copies of $ G $.  Let $ M ^ a := M
    \otimes M \otimes \dots \otimes \ M $ be the tensor product of $ a
    $ copies of $ M $. Then $ M ^ a $ is a $ k G ^ a $-module.
        \par \vspace{5mm}
        Now let $ H := {\rm Sym} ( a ) $. Then $ H $ can be made to
        act on $ G ^ a $ by permuting factors, or as it is sometimes
        called  by {\it place permutations.}    Thus for example we
        have $ ( g_1 , g_2 , \dots , g_a ) \cdot ( 1 \ 2 ) = ( g_2 ,
        g_1 , \dots , g_a ) $.          The action of $ H $ on $ G ^ a
        $ allows us to form     the semidirect product $ G ^ a > \! \!
        \! \! \lhd \ H $ which is known as the wreath product $ G \wr
        H $. The subgroup       $ G ^ a = G \times G \times \dots
        \times G $ is called the base subgroup of $ G $ and the
        subgroup $ H = {\rm Sym} ( a ) $        is called the top
        subgroup of $ G$.
        \par \vspace{5mm}
        Similarly we may define an action of $ H $ on the module $ M ^
        a $ by letting $ H $ permute factors. This turns the    tensor
        product $ M ^ a $ into a $ k H $-module. Therefore $ M ^ a $
        is simultaneously a $ k G ^ a $-module  and a $ k H $-module.
        We may define a $ k ( G \wr H ) $-module structure on $ M ^ a
        $ by stipulating
        $$
        m ^ { g h } = ( m ^ g ) ^ h ,$$
        for $ m \in M ^ a $, $ g
        \in G ^ a $ and $ h \in H $.    We use $ M ^ { \wr a } $ to
        denote this $ k ( G \wr H ) $ module. In essence we let the
        subgroup $ G ^ a $ of $ G \wr H $ act first     via the
        product action and then we let the top subgroup $ H \leq G \wr
        H $ act by permuting factors.
        \par \vspace{5mm}
        The above construction is performed for a $ k G$-module $ M $
        but is equally applicable to the situation of   a permutation
        $ G $-space $ X $. If $ X $ is a permutation $ G $-space then
        we define       a permutation $ G \wr H $-space $ X ^ { \wr a
        } $ as follows. As sets we let $ X ^ { \wr a } $ be the
        cartesian       product $ X \times X \times \dots \times X $
        of $ a $ copies of $ X $. Then we let $ G ^ a $ act via the
        product action  and we let $ H $ act by permuting factors.
        \par 
        If $ M = k X $ is the permutation $ k G $-module afforded by $
        X $ then        $ M ^ { \wr a } = k X ^ { \wr a } $ is the
        permutation module corresponding        to the permutation $ G
        \wr H $-space $ X ^ { \wr a } $. In particular in this
        situation $ M ^ { \wr a } $ is  a permutation $ k ( G \wr H )
        $-module.
        \par \vspace{5mm}
        The construction of the wreath product is of course standard
        and we refer to \cite[(4.1)]{jk81} for more details. Usually
        the     wreath product $ G \wr H $ of two groups $ G $ and $ H
        $ is defined in terms of a permutation $ G $-space $ X $ and
        a permutation $ H $-space $ Y $. In the above construction we
        are taking $ H = {\rm Sym} ( a ) $ and the permutation  $ H
        $-space $ Y $ is the natural permutation $ {\rm Sym} ( a )
        $-space on the set of $ a $ letters $ \{ 1 , 2 , \dots , a \}
        $.  }  {{\bf Tensor Product with Wreath Product Action} }
\end {block}
\mylabel{lemma:DirectProductProjectiveSummand}
\begin {lemma}
  {     Let $ M $ be a $ k G $-module and let $ a $ be a positive
    integer.  The $ k ( G \wr H ) $-module $ M ^ { \wr a } $ admits a
    projective summand only if $ M $ admits a projective summand.  }
  {     We have $ {\rm Res} ^ { G \wr H } _ { G ^ a } ( M ^ { \wr a }
    ) = M ^ a $, so the $ k G ^ a $-module $ M ^ a $ must admit a
    projective summand.  The result now follows.  }
\end {lemma}
\mylabel{lemma:WrProdHasProj}
\begin {lemma}
  {             Let $ M $ be a $ k G $-module and let $ a < p $ be a
    positive integer. Suppose $ M $ admits a projective summand. Then
    the $ k ( G \wr H ) $-module $ M ^ { \wr a} $ admits a projective
    summand.
                
  } {   The size of $ H = {\rm Sym} ( a ) $ is coprime to $ p $
    because $ a < p $, so the base group $ G ^ a $ contains a Sylow $
    p $-subgroup of $ G \wr H $.  We know that the $ k G ^ a $-module
    $ M ^ a $ admits a projective summand.  Now $ M ^ a = {\rm Res} ^
    { G \wr H } _ { G ^ a } ( M ^ { \wr a } ) $ and we may apply
    \arabicc{lemma:ProjectiveModules}{}
    to infer that $ M ^ { \wr a } $ admits a projective summand.
  }
\end {lemma}
\markright{Wreath Products and Projective Summands}
\mylabel{lemma:restrictedproductprojective}
\begin{lemma}
{
    Let $ M $ be a $ k G $-module.
    Assume $ a \geq 2 $  and 
    $ M ^ { \wr a } $ admits a projective
    summand as a $ k G ^ { \wr a } $-module.
    Then $ M ^ { \wr ( a - 1 ) } $ admits a projective
    summand as a $ k G ^ { \wr ( a - 1 ) } $-module.
}
{
    For convenience let $ M_i \cong M $ and
    $ G_i \cong G $ for $ i = 1, \dots a $.
    Then
    $ G ^ { \wr a } \cong ( G_1 \times \dots \times G_{ ( a - 1 ) } \times G_a ) \cdot T $ 
    where $ T \cong {\rm Sym} ( a ) $. So $ G $ contains the
    subgroup
    $ G ^ { \wr ( a - 1 ) } \cong ( G_1 \times \dots \times G_{ ( a - 1 ) } \times 1 ) \cdot T_1 $
    where $ T_1 \leq T $ is the subgroup of $ T $ fixing
    the $ a $-th factor of $ M_1 \otimes \dots \otimes M_a $ 
    and $ T_1 \cong {\rm Sym} ( ( a - 1 ) ) $ .

    The restriction of $ M ^ { \wr a } $ to this subgroup must also
    admit a projective summand and we have
    $$ {\rm Res } ^ { G ^ { \wr a } } _ { G ^ { \wr ( a - 1 ) } } ( M ^ { \wr a } )
       \cong
       \bigoplus _ { { \dim } | M | }  M ^ { \wr ( a - 1 ) },$$
    a direct sum of $ { \dim } | M | $ isomorphic copies
    of $ M ^ { \wr ( a - 1 ) } $.

    It follows that the $ k G ^ { \wr ( a - 1 ) } $-module
    $  M ^ { \wr ( a - 1 ) } $ must have a projective summand.}
\end{lemma}
\setcounter{chapterno}{2} \setcounter{lastplace}{12}
\begin{nblock}
  {     In this section we will define the notion of {\it fixed point
      sets} for a permutation $ G $-space $\Omega $ (over $ k $);
    these will be in a bijective correspondence with the vertices of
    the components of the permutation $ k G$-module $ k \Omega $. The
    vertex corresponding to a given fixed point set is easily
    recovered, and vice versa, and we can reformulate Brou\'e's
    correspondence in terms of this definition. \par \vspace{5mm} The
    maximality property of vertices requires that fixed point sets
    satisfy amongst general subsets of $ \Omega $ a property which we
    will describe as {\it closed}. By restricting our attention to
    closed subsets of $ \Omega $ we may replace in the Brou\'e
    correspondence the quotient group $ N_G ( Q ) / Q $ with another
    quotient group $ N_X / S_X $ which for our purposes is perhaps
    more convenient.
        \par \vspace{5mm}
        The next lemma describes how vertices - and thus projectivity
        - behave when we factor out $ p ^ \prime $-subgroups.    }
      {Main Results on Fixed Point Sets.}
\end{nblock}
\mylabel{lemma:Quotient}
\begin {lemma}
  {     Let $ \Omega $ be a permutation $ G $-space and let $ M := k
    \Omega $ be the corresponding permutation $ k G $-module.  Let $ N
    $ be a normal subgroup of $ G $ that fixes every element of $
    \Omega $ and let $ \pi : G \to G / N $ be the quotient
    homomorphism. Let $ M = M_1 \oplus M_2 \oplus \dots \oplus M_t $
    be a decomposition of $ M $ into indecomposable summands.  Each
    summand $ M_i $ may be viewed as a $ k G $-module or a $ k ( G / N
    ) $ -module.  A subgroup $ Q \leq G $ is a vertex of a summand $
    M_i $ as $ k G $-modules if and only if $ ( Q ) \pi $ is a vertex
    of $ M_i $ when viewed as a $ k (G/N)-module $.  }  {       We may
    assume $ Q \subset P $. Now $ P $ is a Sylow $ p $-subgroup of $
    G$ so $ ( P ) \pi $ is a Sylow $ p $-subgroup of $ G / N $.
    Choose $ i $ and let $ X $ be a $ P $-invariant basis of $ M_i $.
    Then $ X $ is a $ ( P ) \pi $-invariant basis of $ M_i $ viewed as
    a $ k ( G / N ) $-module. If $ Q \leq G $ is another subgroup of $
    G $ then $ Q $ and $ ( Q ) \pi $ fix the same elements of $ X $,
    so the subgroup $ Q $ of $ P $ is maximal subject to fixing an
    element of $ X $ if and only if the subgroup $ ( Q ) \pi $ of $ (
    P ) \pi $ is maximal subject to fixing an element of $ X $. This
    concludes the proof.  }
\end {lemma}
\mylabel{lemma:vertexintersection}
\begin {lemma}
  {     Let $ G $ be a finite group and the finite set $ \Omega $ be a
    transitive permutation $ G $-space. Let $ Q $ be a vertex of some
    indecomposable summand of $ K \Omega $ and let $ P $ be a Sylow $
    p $-subgroup of $ G $ that contains $ Q $. Then there is some
    element $ \omega \in \Omega $ such that $ Q = P \cap \Stab_G (
    \omega ) $.  }  { Let $ k \Omega = M \oplus M ^ \prime $ where $ M
    $ is indecomposable with vertex $ Q \leq P $.  We have $ { \rm
      Res} ^G_P ( k \Omega ) = \bigoplus ^t_{ i = 1 } k \Omega_i $
    where the $ \Omega_i $ are the $ P $-orbits in $ \Omega $.  Note
    that the $ k \Omega_i $ are indecomposable as $ k P $-modules
    since $ P $ is a $ p $-group. By Krull-Schmidt we have therefore
    that $${ \rm Res} ^G_P ( M ) = \bigoplus ^t_{ i \in I } k \Omega_i
    $$
    where $ I $ is some subset of $ \{1, 2, \dots, t \} $. Identify
    $ M $ as a $ P $-module with $ \bigoplus ^t_{ i \in I } k \Omega_i
    $.  Then $ M $ has $ P $-invariant basis $ X = \bigcup_{i \in I}
    \Omega_i $. Let $ Q \leq P $ be a vertex of $ M $. Then by Theorem
    1.2, $ Q $ is a maximal subgroup of $ P $ subject to fixing an
    element of $ X $. So there is some $ \omega \in \Omega_i $ for
    some $ i $ which is fixed by $ Q $ but not by any subgroup  $ R $ of $ P $ properly containing $ Q $. \\
    Then $ P \cap {\rm Stab}_G(\omega) = {\rm Stab}_P(\omega) = Q $.         }
 \end  {lemma}
 \markright{Vertices of $ p $-Permutation Modules}
 \mylabel{theorem:VertIsStab}
\begin {theorm}
  {     Let $ M $ be a permutation $ k G $-module, say $ M = k \Omega $.
    Let $ Q$ be a
    vertex of an indecomposable summand of $ M $, let $ X $ be the set of fixed
    points of $ Q $ on $ \Omega $, and let $ S $ be the pointwise
    stabilizer of $ X $ in $ G $.
    Then $ Q $ is a Sylow $ p $-subgroup of $ S $.  }  {        Let $
    R \geq Q $ be a Sylow $ p $-subgroup of $ S $. If $ R \neq Q $ let
    $ U := N_R ( Q ) $. Then $ U $ is a $ p $-subgroup which strictly
    contains $ Q $ and we have
    $$
    Q <_{ \! \tiny{\neq} } U \leq N_G ( Q ) \cap S .$$
    The factor
    group $ N_G ( Q ) / Q $ therefore contains the non-trivial $ p
    $-subgroup $ \overline{U} := U / Q $. But $ U \leq S $ so $ U $
    and $ \overline{U} $ act trivially on $ F $. Certainly then a
    Sylow $ p $-subgroup of $ N_G ( Q ) / Q $ cannot act regularly on
    any element of $ F $, and by \linebreak Lemma
    \arabicc{lemma:RegularProjective}{} the $ k N_G ( Q ) / Q $-module
    $ k F $ does not admit a projective summand. Of course by
    the Brou\'e correspondence \arabicc{BCorr}{}
    this implies that $ Q $ is not a vertex of an
    indecomposable summand of $ M $. The alternative $ Q = R $ must be
    the case.  }
\end {theorm}
\mylabel{definition:QSM}
\begin {block}
  { Suppose $ \Omega $ is a $ G $-set. Then if $ X \subset \Omega $ is
    a subset we define $S_X := \Stab_G(X) $ to be the pointwise
    stabilizer, and we write $ Q_X $ for a Sylow $ p $-subgroup of $
    S_X $. We will also write $ N_X $ for the set stabilizer of $ X $
    in $ G $, that is, $ N_X := \{ g \in G | X^g = X \} $.  }  { {\bf Notation}}
\end{block}
\begin {definition}
  {      Suppose $ \Omega $ is a $ G $-set and $ X \subseteq Q $. We
    say that $ X $ is a fixed point set if there is a vertex of some
    indecomposable summand of $ k \Omega $ with vertex $ Q $ such that $ X = \Fix_\Omega
    ( Q ) $. \newline Theorem 2.3 shows that in this case $ Q $ is a
    Sylow $ p $-subgroup of $ S_X $ and so for a fixed point set
    we may take $ Q_X = Q $.
    \newline We can reformulate part of the Brou\'e corrspondence as
    follows: } { }
\end {definition}
\begin{corollary}
  { The set $ X \subseteq \Omega $ is a fixed point set if and only if
    the $ p $-permutation module $ k X $ for $ N_G(Q_X)/Q_X $ has a
    projective summand.  }  { }
\end{corollary}

\begin {definition}
  {      Suppose $ \Omega $ is a $ G $-set, and $ X \subseteq \Omega $
    is a subset. We say that $ X $ is {\it closed} if $ \Fix_\Omega (
    Q_X ) = X $. This notion is independent of the choice of $ Q_X $.
    Note that if $ X $ is a fixed point set, then $ X $
    is closed.
        \par 
        One way of obtaining closed subsets is to take the closure of
        a given subset. For subgroups $ L $ of $ G $ we denote by $ *
        $ an assignment (depending on a choice of $ P $) $ L \mapsto
        \Fix_{ \Omega } ( \Syl_p ( L ) ) $ then we mean that    $ L ^
        * $ is the set of fixed points in $ \Omega $ of some Sylow $ p
        $-subgroup of $ L $. The reverse assignment is to denote
        by $ * $ a choice of map $ X \mapsto \Syl_p ( \Stab_G ( X ) )
        $ which assigns to the subset $ X $ of $ \Omega $ some Sylow $
        p $-subgroup    $ X ^ * $ of the stabilizer in $ G $ of $ X $.
        We may say that the subset $ X ^ { * * } \subset \Omega $ is
        the {\it closure} of $ X $.     Thus $ X $ is closed exactly
        when $ X = X ^ { * * } $, and this must happen when for
        instance $ X $ is a fixed point set. Similarly for      a
        vertex $ Q $ of some component of $ k \Omega $ we must have $
        Q = Q ^ { * * } $.  }  { }
\end {definition}
\mylabel{block:PropertiesClosed}
\begin {block}
  { Suppose $ X $ is closed, that is $ \Fix_\Omega ( Q_X ) = X $. Then
    we even have $ \Fix_\Omega (S_X) = X$.  In general, the normalizer
    of a subgroup permutes the fixed points of this subgroup, so the
    normalizer $ N_G(S_X) $ also leaves $ X $ invariant. We have defined
    (in \arabicc{definition:QSM}{})
    $$
    N_X := \{ g \in G : (X)g \subseteq X \} $$
    the largest subgroup
    of $ G $ that leaves $ X $ invariant. In general,
    $ N_X \subseteq N_G( S_X ) $. When
    $X $ is closed, we therefore have $ N_G(S_X) = N_X $.  Our aim now
    is to characterize fixed point sets in terms of $ N_X $ and $ S_X
    $ (which will allow us to characterize vertices in terms of fixed
    point sets).  }  {{\bf Properties of Closed Subsets} }
\end {block}
\markright{Fixed Point Sets}
\begin {block}
  {     Let $ \Omega $ be a $ G $-set and let $ X \subset \Omega $ be
    a closed subset. Define
    $$
    M_X := N_X / S_X .$$
    As we mentioned above, the group $ N_X $
    acts on $ X $ with kernel $ S_X $.  Thus $ M_X $ acts on $ X $ and
    $ X $ is a faithful permutation $ M_X $ -space.  The action of $
    N_X $ on $ X $ induces a homomorphism $ N_X \to {\rm Sym} ( X ) $,
    and $ M_X $ may be taken to be the image of this homomorphism.
  }      {{\bf Notation}         }
\end {block}
\begin {block}
  {      {\it Given a closed subset $ X \subset \Omega $ we have $ N_G
      ( Q_X ) \leq N_X $, and also $ N_G ( Q_X ) \cdot S_X = N_X $.}
    \newline {\bf Proof.} The subgroup $ N_G ( Q_X ) $ will leave
    invariant the set $ \Fix_{ \Omega } ( Q_X ) $; but $ X $ is closed
    so this set is again $ X $, whence $ N_G ( Q_X ) $ leaves $ X $
    invariant giving $ N_G ( Q_X ) \leq N_X $.
         \par 
         That $ N_G ( Q_X ) \cdot S_X = N_X $ follows from a Frattini
         argument since $ Q $ is a Sylow $ p $-subgroup of $ S_X $ and
         $ N_G ( Q_X ) \leq N_X = N_G ( S_X ) $. ( The argument
         attributed to Frattini is that if $ P $ is a Sylow $ p
         $-subgroup of a normal subgroup $ N $ of $ G $ then $ G = N_G
         ( P ) \cdot N $.)  }  {{ \bf Proposition}}
 \end    {block}
 \begin    {theorm}
   {    Let $ X \subset \Omega $ be closed. Then the projective
     components of the $ k M_X $-module $ k X $ are in bijective
     correspondence with the projective components of the $ k ( N_{ G
     } ( Q_X ) / Q_X ) $-module $ k X $.  }  {  By the preceding
     proposition we have
     $$
     M_X = N_X / S_X = N_G ( Q_X ) \cdot S_X / S_X $$
     $$
     \cong N_G ( Q_X ) / N_{ S_X } ( Q_X ) \cong ( N_{ G } ( Q_X )
     / Q_X ) / ( N_{ S_X } ( Q_X ) / Q_X ) .$$  Now $ N_{ S_X } ( Q_X
     ) / Q_X $ has $ p ^ \prime $ size since $ Q_X $ is a Sylow $ p
     $-subgroup of $ S_X $ . The result follows from Lemma \arabicc
     {lemma:Quotient}{}.  }
\end {theorm}
\vspace*{0.1cm}
\par We can now reformulate the Brou\'e Correspondence:
\begin {block}
  {     {\bf Brou\' e's correspondence for closed sets} Let $ X
    \subset \Omega $ be closed. There is a multiplicity-preserving
    bijective correspondence between the isomorphism types of
    components of $ k \Omega $ which admit $ X $ as a fixed point set
    and the set of projective components of the permutation $ k M_X
    $-module $ k X$ (up to isomorphism).  }  {  {\bf Theorem.} }
\end {block}
\mylabel{corollary:FixIffProjComp}
\begin {block}
  {     {\it Let $ X \subset \Omega $ be closed. Then $ X $ is a fixed
    point set if and only if the permutation $ k M_X $-module $ k X $
    admits a projective component. } }  {        {\bf Corollary.}  }
\end {block}
 
\setcounter{chapterno}{3}
 
\newcommand{\SetS}{ {\mathcal SetS } } \newcommand{\starbar}{ * ^ { \!
    \! \! \! \! \! } }

\setcounter{lastplace}{1}
\markboth{Building Fixed Point Sets}{Building Fixed Point Sets}
\begin {nblock}
  {      {\it A Specific Problem.}
    We have seen that the problem of
    determining the vertices of the indecomposable summands of
    permutation modules of finite groups is equivalent to finding the
    corresponding fixed point sets. Now we shall study these fixed
    point sets for the permutation $ k {\rm Sym} ( q m ) $-module $ k
    \Xi^q_{qm} $ where $ \Xi^q_{qm} $ is the conjugacy class of $ {\rm
      Sym} ( q m ) $ consisting of fixed point free elements which
    are products of $ q $-cycles and $ k $ is a field of
    characteristic $ p $.  }  {{\rm Sym}metric Groups and Conjugacy
    Actions.}
\end {nblock}
\begin {block}
  {     In such actions, that is when a group $ G $ acts on itself by
    conjugation, the subsets of the relevant $ G $-set have an
    additional structure given by the group multiplication. We will
    therefore consider ways of building new fixed point sets out of
    known ones, using the provided multiplication, and also using the
    fact that the direct product of two symmetric groups can be
    embedded into a larger symmetric group.
            \par In the first part of this section we will study set constructions, and at the end of the section we will relate these
            set constructions to fixed point sets.  }  {{\bf
              Irreducible Sets} }
\end {block}
\mylabel{definition:FS}
\begin{definition}
  {     The group of finitary permutations of $ {\mathbb N } $ is
    defined to be the union of all finite symmetric groups:
    $$
    \Sym_{f} ( {\mathbb N }) := \bigcup_{n=1}^{\infty} \Sym ( n ) .$$
    Recall that the support of a permutation $ \rho $ is defined as:
    $$
    \supp \ \rho := \{ \omega \in {\mathbb N} \ | \ \omega ( \rho )
    \neq \omega \} .$$
    If $ X $ is a set of finitary permutations we
    set $ supp(X) = \bigcup_{\rho \in X} \supp (\rho) $. Now we define
    $ {\mathcal F} $ to be all finite subsets of $ \Sym_f ( {\mathbb N
    }) $, excluding the empty set and the set which contains only the
    identity permutation.  That is,
    $$
    {\mathcal F} := \{ X \subset \Sym_f ( {\mathbb N }) : 0 < | X |
    < \infty \ {\rm and} \ X \neq \{ 1 \} \} .$$
            
    For an integer $ q \geq 2 $, we define $ {\mathcal S} ^ q
    \subset {\mathcal F} $ to be all subsets $ X $ such that  \vspace{0.5cm} \newline    
    \hspace*{10 mm} (1) each $ \rho \in X $ is a product of $ q $-cycles; \newline
    \hspace*{10 mm} (2) for each $ \rho \in X $ we have $ \supp ( \rho ) = \supp ( X )
    $. \newline
           
    Any such $ X \in {\mathcal S} ^ q $ is conjugate to a subset of $
    \Xi_{qm}^q $ where $ qm $ is the size of $ \supp ( X ) $.  }  { }
\end {definition}
\mylabel{notation:QSMagain}
\begin {block}
  {     Let $ X $ be an element of $ {\mathcal F} $. By construction $
    X $ is a subset of some finite subgroup of $ {\rm Sym} ( {\mathbb
      N} ) $, and in fact we have $ X \subset {\rm Sym} ( \supp \ X )
    \leq {\rm Sym} ( {\mathbb N} ) $. So let $ G_X := {\rm Sym} ( supp
    \ X ) $, and write $ N_X := N_{ G_X } ( X ) = \{ g \in G_X \ | X ^
    g = X \} $ and $ S_X := \Stab_{ G_X } ( X ) $. \par Now $ S_X $ is
    a normal subgroup of $ N_X $ and $ N_X $ acts by conjugation on $
    X $ with kernel $ S_X $. Thus $ M_X := N_X / S_X $ acts on $ X $
    with trivial kernel and $ k X $ is a faithful permutation $ k M_X
    $-module. Let $ Q_X \leq S_X $ be a Sylow $ p $-subgroup of $ S_X
    $.
        \par 
        These definitions were also made in Section
        \arabicc{definition:QSM}{}, and it was suggested there that a
        distinctive case is when the    permutation $ k M_X $-module $
        k X $ admits a projective summand. This is the condition which
        characterizes, amongst closed   sets, the fixed point sets.  }
      {{\bf Notation} }
\end {block}
\mylabel{block:MultDegree}
\begin{block}
  {      Let $ X, Y \in {\mathcal F}.$ Letting $ a := | \supp \ X | $
    and $ b := | \supp \ Y | $, there is an embedding of groups
    $$
    \tau : {\rm Sym} ( a ) \times {\rm Sym} ( b ) \to {\rm Sym} ( a
    + b ). $$  We identify $ {\rm Sym} ( a ) $ with $ {\rm Sym} ( supp
    \ X ) $ and $ {\rm Sym} ( b ) $ with $ {\rm Sym} ( \supp \ Y ) $ so
    that this embedding restricts to an embedding of the Cartesian
    product of $ X $ and $ Y $:
    $$
    \tau : X \times Y \to {\rm Sym} ( a + b ) .$$
    We write $ X * Y
    $ to mean the image of this embedding. By convention we write $ X
    * \{ 1 \} = X $ for any subset $ X \in {\mathcal F} $.  Now
    consider the $ s $-fold Cartesian product $ X ^ s = X \times X
    \times \dots \times X $ and the embedding $ \tau : X ^ s \to {\rm
      Sym}( a s ) $. We let $ \Delta ^ s X $ be the image under $ \tau
    $ of the diagonal subset $ \{ ( x , x , \dots , x ) \ | \ x \in X
    \} $ of $ ( X ) ^ s $.
        \par 
        We write the {\bf degree} of $ Y $ as $ d ( Y ) $ and define
        it by $ d ( Y ) := | \supp \ Y | $.  }  {{\bf Multiplication In
          $ {\mathcal F} $} }
\end {block}
\markright{Definition of $ * $-Product}
\begin {remark}
  {     We will say that two elements $ X $ and $ Y $ of $ {\mathcal
      F} $ are equivalent if there is $ \sigma \in {\rm Sym} (
    {\mathbb N} ) $ such that $ X ^ \sigma = Y $. The $ * $-product
    can be thought of as defined on the resulting equivalence classes
    and in this case the $ * $-product is abelian and associative.  }
  { }
\end {remark}
\begin {definition}
  {     {\bf Irreducible Sets.} We say that $ X \in {\mathcal F} $ is
    {\bf irreducible} if there do not exist elements $ Y , Z \in
    {\mathcal F} $ such that $ X = Y * Z $. Clearly $ X $ can be
    written as a product of irreducible members of $ {\mathcal F} $.
    Note that $ X \in {\mathcal S} ^ q $ is irreducible if and only if
    there do not exist in $ Y, Z \in {\mathcal S} ^ q $ such that $ X
    = Y * Z $.  }  { }
\end {definition}
\begin {remark}
  {     If we interpret the product $ * $ as being defined on the set
    of equivalence classes of $ {\mathcal F} $ then Lemma 3.9
    immediately implies that any decomposition into irreducible
    factors is unique.  In fact slightly more is implied because the
    statement "$ Y_i = Z_i $" which forms part of the statement of the
    lemma asserts an absolute equality, that is, as subsets of $ {\rm
      Sym} ( {\mathbb N} ) $, and not just equivalence. }  { }
\end {remark}

\mylabel{lemma:FixCombo}
\begin {lemma}
  {     The set $ Y \in {\mathcal S } ^ q $ is a fixed point set if
    and only if each of the following hold
        \begin {enumerate}
        \item $ Y $ is closed
        \item the permutation $ k M_Y $-module $ k Y $ admits a
          projective summand.
        \end {enumerate}
      } {       This is a combination of Corollary
        \arabicc{corollary:FixIffProjComp}{} and Theorem
        \arabicc{theorem:VertIsStab}{} .  }
\end {lemma}
\markright{Irreducible Fixed Point Sets}
\mylabel{lemma:UniqueDecomp}
\begin {lemma}
  {     {\bf Uniqueness of Irreducible Decomposition.} Let $ X \in
    {\mathcal F} $ be such that $ X = Y_1 \times Y_2 \times \dots
    \times Y_s = Z_1 \times Z_2 \times \dots \times Z_t $ where the $
    Y_i $ and $ Z_i $ are irreducible elements of $ {\mathcal F} $.
    Then $ s = t $ and after a reordering of factors we have $ Y_i =
    Z_i $ for $ 1 \leq i \leq s $.  }  {        We may assume that
    respectively the $ Y_i $ and the $ Z_i $ have pairwise disjoint
    supports and that $ X = Y_1 \times Y_2 \times \dots \times Y_s =
    Z_1 \times Z_2 \times \dots \times Z_t $. Put $ s_i := \supp \ Y_i
    $ and $ t_i := \supp \ Z_i $. Now $ X $ is a subset of the subgroup
    $ {\rm Sym} ( s_1 ) \times {\rm Sym} ( s_2 ) \times \dots \times
    {\rm Sym} ( s_s ) \leq {\rm Sym} ( {\mathbb N} ) $ and we may
    consider the projection $ \pi_ { Y_i } : X \to Y_i \subset {\rm
      Sym} ( s_i ) $, and similarly the projection $ \pi_ { Z_i } : X
    \to Z_i $. Thus we have $ Y_1 = ( X ) \pi_ { Y_1 } = ( Z_1 \times
    Z_2 \times \dots \times Z_t ) \pi_ { Y_1 } = ( Z_1 ) \pi_ { Y_1 }
    \times ( Z_2 ) \pi_ { Y_1 } \times \dots \times ( Z_t ) \pi_ { Y_1
    } $ and since $ Y_i $ is irreducible this forces $ Y_1 = ( Z_i )
    \pi_ { Y_1 } $      for some $ i $, and $ s_1 \subset t_i $.
         \par 
         Similarly we have $ Z_i = ( Y_1 ) \pi_ { Z_i } \times ( Y_2 )
         \pi_ { Z_i } \times \dots \times ( Y_s ) \pi_ { Z_i } $. But
         $ s_1 \subset t_i $ so $ ( Y_1 ) \pi_ { Z_i } = Y_1 $, and
         the irreducibility of $ Z_i $ implies $ Z_i = Y_1 $.
         Reordering the $ Z_j $ we have $ Z_1 = Y_1 $; upon cancelling
         these two factors the result follows by induction.  }
\end {lemma}
\vspace*{0.1cm}
\par Let $ N_Y = \{ g \in {\rm Sym} (\supp \ Y) : Y ^ g = Y \} $  for $ Y \in {\mathcal F} $.
\mylabel{corollary:PermutesFactors}
\begin {corollary}
  {     Let $ Y \in {\mathcal F} $ and let $ Y = Y_1 \times Y_2 \times
    \dots \times Y_s $ be a decomposition into irreducibles. For every
    $ 1 \leq i \leq s $ and $ g \in N_Y $ there is $ 1 \leq j \leq s $
    such that $ Y_i ^ g = Y_j $.  }  {  We have $ Y = Y ^ g $ and $
    Y_1 \times Y_2 \times \dots \times Y_s = Y_1 ^ g \times Y_2 ^ g
    \times \dots \times Y_s ^ g $ so the result follows immediately
    from the lemma.  }
\end {corollary}
\begin {corollary}
  {     Then the action of $ N_Y $ on $ Y $ induces a permutation $
    N_Y $-action on the set of irreducible factors of $ Y $.  }
  {     This is a restatement of the previous corollary.  }
\end {corollary}
\mylabel{lemma:DeltaIsIrreducible}
\begin  {lemma}
  {     Suppose that $ X \in {\mathcal S } ^ q $ is irreducible and
    satisfies $ | X | > 1 $.  Then $ \Delta ^ i X $ is irreducible for
    every $ i $.  }  {  We write $ \Delta ^ i X = X_1 \Delta X_2
    \Delta \dots \Delta X_i $ where each $ X_j $ is conjugate to $ X $
    and where the $ \supp \ X_j $ are pairwise disjoint. Put $ \alpha
    := \bigcup_{ j = 1} ^ i \supp \ X_j $.  Assume that $ \Delta ^ i X
    = Y \times Z $ where $ Y $ and $ Z $ are elements of $ { \mathcal
      S ^ q} $ with disjoint non-empty support. Then for $ 1 \leq j
    \leq i $ we have $ X_j = ( \Delta^i X ) \pi_{ X_j } = ( Y ) \pi_{
      X_j } \times ( Z ) \pi_{ X_j } $, and since $ X_j $ is
    irreducible this implies that either $ X_j = ( Y ) \pi_{ X_j } $
    or $ X_j = ( Z ) \pi_{ X_j } $.  Consequently either $ \supp \ X_j
    \subset \supp \ Y $ or $ \supp \ X_j \subset \supp \ Z $.
        \par 
        Let $ a $ be the set of         $ j $ such that $ \supp \ X_j
        \subset \supp \ Y $ and let $ b $ be the set of $ j $ such that
        $ \supp \ X_j \subset \supp \ Z $. Since  $ \supp \ Y $ and $
        \supp \ Z $ are non-empty it follows that $ a $ and $ b $ are
        each non-empty. Also, since     $ \alpha = \supp \ Y \cup supp
        \ Z $ we have $ \supp \ Y = \bigcup_{ j \in a } \supp \ X_j $
        and $ \supp \ Z = \bigcup_{ j \in b } \supp \ X_j $.      After
        a possible reordering we may assume that $ a = \{ 1 , 2 ,\dots
        , s \} $ and $ b = \{ s + 1 , s +2 , \dots , i \} $.
        \par 
        Comparing this with $ X_1 \Delta X_2 \Delta \dots \Delta X_i =
        Y \times Z $ we have    $ Y = ( Y \times Z ) \pi_Y = X_1
        \Delta X_2 \Delta \dots \Delta X_s $ and similarly      $ Z =
        X_{ s + 1} \Delta X_{ s + 2 } \Delta \dots \Delta X_i $. Thus
        $ | Y | = | Z | = | X | $ and $ | Y \times Z | = | X | ^2 $.
        But $ | Y \times Z | = | \Delta ^ i X | = | X | $, and so $ |
        X | = 1 $.  }
\end {lemma}
\begin{block}
  { Let $ X \in {\mathcal F} $, and define
    $$
    G_X := \Sym ( \supp \ X) ,$$
    the symmetric group permuting the
    elements in $ \supp \ ( X ) $.
     \par 
     We view $ X $ as a subset of a permutation $ G $-set, where $ G $
     acts by conjugation. That is,
     $$
     \Omega = \Xi_X := \bigcup_{ g \in G } X ^ g \subseteq G_X . $$
     We observe that in this setup $ S_X $ is equal to $ C_{ G_X } ( X ) $
     and $ N_X $ is equal to $ N_{ G_X } ( X ) $. The set of fixed points
     of $ Q_X $ on $ \Omega $ is then equal to
     $$ {\rm Fix}_{ \Omega } ( Q_X ) = C_{ G_X } ( X ) \cap \Omega .$$
     Note that we always take these with respect
     to the $ G_X $-set just defined.  }  {{\bf Conjugation actions.}  }
\end{block}
\mylabel{block:FPS}
\begin {block}
  {      Now we specialize these concepts to $ Y \in {\mathcal S} ^ q
    $. Recall then that $ Y $ is a fixed point set if there is an
    indecomposable summand of the permutation $ k G_Y $-module $ k
    \Xi_Y $ with vertex $ Q \leq G_Y $ such that $ Y = \Fix_{\Xi_Y} (
    Q ) $. Thus $ Y $ is a fixed point set if and only if $ Q_Y $ is a
    vertex of an indecomposable summand of the permutation module $ k
    \Xi_Y $, for the group $ G_Y $. For example, if $ Y \in {\mathcal S} ^ q $
    is such that all elements of $ Y $ have the same support then $ \Xi_Y $
    is actually the $ G_Y $ set of all fixed point free elements which
    are the product of $ q $-cycles.
    
    }  {{\bf Conventions} }
\end {block}
\mylabel{block:Exact}
\begin {block}
  {     Let $ Y \in {\mathcal S } ^ q $. We say that $ Y $ is {\bf
      exact} if $ Q_Y $ is fixed point free in its action on $ \supp \ 
    Y $, that is, if $ \supp \ Q_Y = \supp \ Y $. Recall that
    for $ Y, Z \in {\mathcal S } ^ q $ we have make choices
    so that $ Q_X *Q_Y \leq Q_{ X * Y } $, and it follows that
    the product of exact sets is exact.
    \par On the other hand, we
    say that $ Y $ is {\bf projective} if $ Q_Y = 1 $. We observe that
    $ Y \in {\mathcal S} ^ q $ is projective and closed if and only if
    $ Y =  \Xi_Y $.
   }  {{\bf Exact Sets
          and Projective Sets} }
\end {block}
\mylabel{lemma:ExactOrProjective}
\begin {lemma}
  {     Suppose that $ Y \in {\mathcal S } ^ q $ is irreducible and
    closed. Then $ X $ is exact or $ X $ is projective.  }  {   Let $ \alpha
    = \supp \ Q_X $ and write $ \supp ( X ) = \alpha \cup  \beta $
    where $ Q_X $ fixes each element of $ \beta $.
    The group $ \langle X \rangle $ centralizes $ Q_X $ and therefore
    it permutes the orbits of $ Q_X $ on the support of $ X $.
    Assume for a contradiction that $ \alpha $ and $ \beta $ are both
    $ \neq \emptyset $.
    Then $ X \subset \Sym ( \alpha ) \times \Sym ( \beta ) $.
    \par
    So we have $ X \subseteq X_\alpha * X_\beta $
    (where $ X_\alpha $ is the set of all factors of $ X $ which
    are supported on $ \alpha $.
    \par
    By assumption, $ X $ is closed, that is
    $ X = \Fix_{ \Omega_X } ( Q_X ) $.
    \par
    We can describe these fixed points: they are all elements
    in $ \Omega_X $ which lie in
    $ \Sym ( \alpha ) \times \Sym ( \beta ) $, and where the
    factor in $ \Sym ( \alpha ) $ commutes with $ Q_X $.
    The factor in $ \Sym ( \beta ) $ is arbitrary. In fact we
    get $ X = Y * X_\beta $ where
    $$  Y = \Fix_{ \Omega_{ X_\alpha } } ( Q_X ), X_\beta = \Omega_{ X_\beta } .$$
    We assume $ X $ is irreducible, so this is a contradiction.
      }
\end {lemma}
 \vspace*{2mm} \newline {\bf Remark. } { We have seen at the same time that if $ Q_X = 1 $
  then $ X = X_\beta $ and this is the complete $ G_X $-set $ \Omega_X $.
  This will be relevant later. }
\setcounter{chapterno}{4} \setcounter{lastplace}{15}
\begin{nblock}
  {      Suppose that we have found a fixed point set $ Y \in
    {\mathcal S } ^ q $ which is irreducible and exact.  By Lemma
    \arabicc{lemma:FixCombo}{} $ Y $ is closed and the permutation $ k
    M_Y $-module $ k Y $ admits a projective summand.  We wish to find
    conditions under which the elements $ * ^ s Y $ and $ \Delta ^ s Y
    $ are fixed point sets. This means that (i) we want find
    conditions under which $ * ^ s Y $ or $ \Delta ^ s Y $ is a closed
    element of $ {\mathcal S } ^ q $, and (ii) we would like to know
    the structure of the permutation $ k M_{ * ^ s Y } $-module $ k *
    ^ s Y $ and of the permutation $ k M_{ \Delta ^ s Y } $-module $ k
    \Delta ^ s Y $. We address the second of these first. In doing so
    we shall make use of the wreath product constructions described in
    Section \arabicc{section:WreathProd}{}.  }  { Permutation
    Structure of Products of Irreducible Sets.  }
\end {nblock}
\mylabel{block:ActOnProd}
\begin {block}
  {     It is well known that the stabilizer in $ {\rm Sym} ( q n ) $
    of an element of $ \Xi_{qn} ^ q $ (that is, a fixed point free
    element which is a peoduct of $ q $-cycles) is isomorphic to the wreath product $ {
      \mathbb Z}_q \wr {\rm Sym} ( n ) $.  We now consider the version
    of this result which applies to the situation where we are dealing
    not with $ q $-cycles in $ {\rm Sym} ( q n ) $ but equal sized
    subsets of $ \{ 1 , 2 , \dots , q n \} $, or of some similar set.
         \par 
         Let $ \alpha $ be a finite subset of $ {\mathbb N} $ and let
         the subsets $ \alpha_1 , \alpha_2, \dots , \alpha_s $
         partition $ \alpha $ and satisfy $ | \alpha_i | = | \alpha_1
         | $. Let $ G := {\rm Sym} ( \alpha_1 ) $. Then $ \alpha_1 $
         is a natural permutation $ G $-space and as described in
         Section \arabicc{section:WreathProd}{} we may form the tensor
         product $ \alpha_1 ^ { \wr s} $ giving a permutation
         $ ( G \wr {\rm Sym} ( s ) ) $-space.
         In fact this construction can
         be made inside $ {\rm Sym} ( \alpha ) $. The action of
         $ {\rm Sym} ( \alpha ) $ on $ \alpha $
         induces an action on sets
         of subsets of $ \alpha $ and we may define $ C $ to be the
         stabilizer under this action of the set $ A := \{ \alpha_1 ,
         \alpha_2, \dots , \alpha_s \} $. It is reasonably clear that
         $ C $ is isomorphic to the wreath product $ {\rm Sym} (
         \alpha_1 ) \wr {\rm Sym} ( s ) $. However, we shall now
         present the details of the proof.  }  {{\bf Actions on
           Products} }
\end {block}
\mylabel{section:FurtherWrProd}
\begin {block}
  {     The base group in the situation above is the Young subgroup $
    B := {\rm Sym} ( \alpha_1 ) \times {\rm Sym} ( \alpha_2 ) \times
    \dots \times {\rm Sym} ( \alpha_s ) $ and there is a canonical
    choice of subgroup $ T \leq C , T \cong {\rm Sym} ( s ) $ whose
    action on the set $ A $ is the full {\rm Sym}metric action on $ s
    $ elements. This canonical choice is determined by the following.
    For each $ i , j $ there is a unique order preserving involution $
    g_{ i j } \in {\rm Sym} ( {\mathbb N} ) $ such that $ \supp \ g_{ i
      , j } = \alpha_i \cup \alpha_j $ and $ \alpha_i ^ { g_{ i j } }
    = \alpha_j $. We have $ g_{ i j } = g_{ j i } $
    and $ g_{ i j } $ can be viewed as the
    transposition $ ( \alpha_i \alpha_j ) $ in $ {\rm Sym} ( A ) $.
    Thus the $ g_{ i j } $ generate a subgroup $ T \cong {\rm Sym} ( A
    ) \cong {\rm Sym} ( s ) $ of $ C $.
        \par 
        Note that although $ T $ is transitive on $ A $,
        each of the above generators of $ T $
        is order-preserving so the set of least elements of the
        $ \alpha_i $ is left invariant by $ T $, and $ T $ is
        not transitive on
        $ \alpha = \alpha_1 \cup \dots \cup \alpha_s $
        if $ | \alpha_1 | > 1 $.

        \par
        The subgroup
        $ B \cdot T = {\rm Sym} ( \alpha_1 ) \cdot T \leq C $ is
        isomorphic to   the wreath product $ {\rm Sym} ( \alpha_1 )
        \wr T $. If $ g \in {\rm Sym} ( \alpha ) $ is any other
        element such that $ A ^ g = A $ then    $ g $ induces some
        permutation on the elements of $ A $ and there is an element $
        t \in T = {\rm Sym} ( A ) $ which is such that  $ g t $ leaves
        each $ \alpha_i $ invariant. Thus $ g t $ lies in the base
        group $ B $ and we have         $ C = B \cdot T \cong {\rm
          Sym} ( \alpha_1 ) \wr {\rm Sym} ( s ) $.        
      }
{{\bf Wreath Products} }
\end {block}
 \begin  {block}
   {    If $ Y $ is an element of $ { \mathcal S ^q } $ then the
     product $ \large { * } ^ s Y $ is equivalent to any product $ Y_1
     \times Y_2 \times \dots \times Y_s $ where each $ Y_i $ is
     conjugate in $ {\rm Sym} ( {\mathbb N} ) $ to $ Y $ and where the
     sets $ \supp \ Y_i $ partition $ \supp \ \large { * } ^ s Y $.  Thus by
     taking $ \alpha _i : = \supp \ Y_i $ and $ \alpha := \supp \ \large
     { * } ^ s Y $ we are in the situation of the preceding paragraphs
     and there is a canonical choice of subgroup $ T \cong {\rm Sym} (
     s ) $ in $ G_{ \large { * } ^ s Y } = {\rm Sym} ( \supp \ \large {
       * } ^ s Y ) $ which permutes the $ \alpha_i $. We may replace $
     Y_i \subset {\rm Sym} ( \alpha_i ) $ with any subset of $ {\rm
       Sym} ( \alpha_i ) $ equivalent to $ Y_i $. We may therefore
     assume that $ Y_i = Y_1 ^ { g_{ 1 i } } $. The subgroup $ T $ now
     permutes the factors $ Y_i $. In the following lemmas we shall
     assume that the subgroups $ T $, $ B $ and $ C $ are as defined
     above (recall also the conventions given in
     \arabicc{block:FPS}{}
     and \arabicc{block:Exact}). }  {{ \bf Wreath Products and Irreducibles} }
 \end  {block}
 \mylabel{lemma:StructureNProd}
\begin {lemma}
  {      Let $ Y \in {\mathcal S } ^ q $ be irreducible. Then
    $$
    N_{ \large{ * } ^ s Y } = ( * ^ s N_Y ) \cdot T \cong N_Y \wr
    {\rm Sym} ( s ) .$$
  } {   Write $ \large{ * } ^ s Y = $ $ Y_1
    \times Y_2 \times \dots \times Y_s $ where each $ Y_i $ is
    conjugate in $ {\rm Sym} ( {\mathbb N} ) $ to $ Y $ and where the
    $ \supp \ Y_i $ are pairwise disjoint.  On one hand we clearly have
    $ T \leq N := N_{ \large{ * } ^ s Y } $., and also $ * ^ s N_Y
    \subseteq N $.
    
    Now let $ g \in N $. The irreducibility condition implies that $ g
    $ must permute the factors $ Y_i $ and so $ g $ must permute the
    subsets $ \supp \ Y_i $. In other words $ g $ must lie in the
    subgroup $ C $. There is an element $ t \in T \leq N $ such that $
    g t $ leaves each $ \supp \ Y_i $ invariant. Thus $ g t \in B \cap
    N = \large{ * } ^ s N_Y $ and $ N \leq ( \large{ * } ^ s N_Y )
    \cdot T $. Since $ N_Y \leq N $ we have $ N = ( * ^ s N_Y ) \cdot
    T \cong N_Y \wr {\rm Sym} ( s ) $.  }
\end {lemma}
\markright{Products of Irreducible Fixed Point Sets}
\mylabel{lemma:StructureS}
\begin {lemma}
  {      Let $ Y \in {\mathcal S } ^ q $ be irreducible. Then
    $$
    S_{ \Delta ^ s Y } = ( * ^ s S_Y ) \cdot T \cong S_Y \wr {\rm
      Sym} ( s ) .$$
    
  } {   Write $ S := S_{ \Delta ^ s Y } $. Then $ T \leq S $ and $ B
    \cap S = \large{ * } ^ s S_Y $ so the proof of the last lemma
    works here.  }
\end {lemma}
\par If $ Y \in {\mathcal S } ^ q $ contains exactly one element then $  \large{ * } ^ s Y = \Delta ^ s Y $. The next lemma assumes that this is not the
case.
 \mylabel{lemma:StructureSProd}
\begin {lemma}
  {     Let $ Y \in {\mathcal S } ^ q $ be irreducible with $ | Y | >
    1 $. Then
    $$
    S_{ \large{ * } ^ s Y } = \large{ * } ^ s S_Y .$$
  } {   We may
    assume $ s \geq 2 $. Write as usual $ \large{ * } ^ s Y = Y_1
    \times Y_2 \times \dots \times Y_s $ and put $ S := S_{ \large{ *
      } ^ s Y } $. Clearly we have $ \large{ * } ^ s S_Y \leq S $.
    Let $ a $ and $ b $ be distinct elements in $ Y $ and let $ g \in
    S $. Then $ \Delta ^ s Y \subset \large { * } ^ s Y $ so $ g $
    certainly fixes $ \Delta ^ s Y $, and $ g \in S_{ \Delta ^ s Y } =
    (* ^ s S_Y) \cdot T $, using the previous lemma.  Thus $ g = u t $
    where $ u \in \large { * } ^ s S_Y \leq S $ and $ t \in T $, forcing $ t \in S $.
    Suppose that $ Y_1 ^ t = Y_i$. If $ i \neq 1 $ then the element $
    a * b * \dots * b \in * ^s Y $, corresponding to the `coordinate
    vector' with an $ a $ in the first place and $ b $'s elsewhere, is
    not centralized by $ t $. It follows that $ Y_1 ^ t = Y_1 $;
    similarly $ Y_i ^ t = Y_i $ for each $ i $ and $ t = 1 $.
    Therefore $ S \leq \large{ * } ^ s S_Y $ as required.  }
\end {lemma}
 \begin  {lemma}
   {    Let $ Y \in { \mathcal S^q} $ be irreducible and let $ S :=
     S_{ \Delta ^ s Y } $. Then
     $$
     N_{ \Delta ^ s Y } = \langle \Delta ^ s N_Y , S \rangle .$$
   }
   {    Put $ N := N_{ \Delta ^ s Y } $. Clearly $ \Delta ^ s N_Y \leq
     N \leq C $ and $ T \leq S \leq N $.
     Let $ g \in N $; then $ g = u
     t $ where $ u \in B $ and $ t \in T $.  Since $ t \in N $ we get $ u \in N $.
     Writing all elements in $
     B = \large { * } ^ s G_Y $ as vectors with coordinates in $ G_Y $
     we have $ u = ( u_1 , u_2 , \dots , u_s ) $. If $ y \in Y $ then
     $ ( y , y , \dots , y ) \in \Delta ^ s Y $ and $ ( y , y , \dots
     , y ) ^ u = ( y ^ { u_1 } , y ^ { u _2 } , \dots , y ^ { u_s } )
     \in \Delta ^ s Y $. Therefore $ y ^ { u_1 } \in Y $ so $ u_1 \in
     N_Y $. Also $ y ^ { u_i } = y ^ { u_1 } $ holds for every $ y \in Y
     $ so $ u_1 u_i ^ { - 1 } \in S_ Y $ and
     $ u = (  u_1 u_1^ { - 1 },  u_2 u_1^ { - 1 } , \dots , u_s u_1 ^ { -1 } ) ( u_1 , u_1 , \dots , u_1 ) $
     gives $ u \in \Delta ^ s N_Y \cdot \large { * } ^ s
     S_Y \leq \Delta ^ s N_Y \cdot S $, as needed.  }
\end {lemma}
\par 
\begin {remark}
  { As a corollary to this lemma we have $ N_{ \Delta ^ s Y } / S_{
      \Delta ^ s Y } $ \linebreak $ \cong \Delta ^ s N_Y / ( \Delta ^
    s N_Y \cap S ) = \Delta ^ s N_Y / \Delta ^ s S_Y \cong N_Y / S_Y
    \cong M_Y $. Thus $ M_{ \Delta ^ s Y } \cong M_Y $. We want to
    show slightly more than this, namely that the pair $ ( M_{ \Delta
      ^ s Y } , \Delta ^ s Y ) $ is permutation isomorphic to the pair
    $ ( M_Y , Y ) $.
        \par 
        If $ Y \in {\mathcal S } ^ q $ we shall let $ \tau : G_Y \to
        \Delta ^ s G_Y $ be the bijection given by      $ g \mapsto (
        g , g , \dots , g ) $.  }  { }
\end {remark}
\mylabel{lemma:StructureOfDeltaProduct}
\begin {lemma}
  {     Let $ Y \in { \mathcal S ^q } $ be irreducible.  Then the pair
    $ ( M_{ \Delta ^ s Y } , \Delta ^ s Y ) $ is permutation
    isomorphic to the pair $ ( M_Y , Y ) $.  }  {       Let $
    \pi_\Delta : N_{ \Delta ^ s Y } \to {\rm Sym} ( \Delta ^ s Y ) $
    be the homomorphism induced by the action of $ N_{ \Delta ^ s Y }
    $ on $ \Delta ^ s Y $ and let $ \pi_Y : N_Y \to \Sym (Y) $ be the homomorphism
    induced by the action of $ N_Y $ on $ Y $. By the previous lemma
    we may let $ \iota : \Delta ^ s N_Y \to N_{ \Delta ^ S Y } $ be
    the inclusion homomorphism.  Then there is a homomorphism
    $$
    \tau \iota \pi_\Delta : N_Y \to \Delta ^ s N_Y \to N_{ \Delta ^
      S Y } \to (N_{ \Delta ^ S Y }) \pi_\Delta = M_{ \Delta ^ S Y }
    \leq {\rm Sym} ( \Delta ^ s Y ), $$
    The homomorphism $ \tau \iota
    \pi_\Delta : N_Y \to M_{ \Delta ^ S Y } $ is surjective since by
    the previous lemma the kernel $ S_{ \Delta ^ S Y } $ of $
    \pi_\Delta $ together with $ \Delta ^ s N_Y $ generates $ N_{
      \Delta ^ S Y } $.
        \par 
        Now $ \tau $ restricts to a bijection $ \tau : Y \to \Delta ^
        s Y $ which induces a isomorphism       $ \tau : {\rm Sym} ( Y
        ) \to {\rm Sym} ( \Delta ^ s Y ) $ and there is a homomorphism
        $$
        \pi_Y \tau : N_Y \to {\rm Sym} ( Y ) \to {\rm Sym} ( \Delta
        ^ s Y ) .$$
        In fact for all $ g \in N_Y $ the elements $ (
        g ) \tau \iota \pi_\Delta $ and $ ( g ) \pi_Y \tau $ are each
        equal to that   element of $ {\rm Sym} ( \Delta ^ s Y ) $ which
        takes $ ( y , y , \dots , y ) \in \Delta ^ s Y $ to $ ( y ^ g
        , y ^ g , \dots , y ^ g ) \in \Delta ^ s Y $.   Thus $ \tau
        \iota \pi = \pi_Y \tau $. But $ M_{ \Delta ^ s Y } = N_{
          \Delta ^ s Y } \pi_{\Delta} = N_Y \tau \iota \pi_{\Delta} $ and $ M_Y = N_Y
        \pi_Y $ so      $ M_Y \tau = M_{ \Delta ^ s Y } $.  }
\end {lemma}
\mylabel{lemma:StructureOfProduct}
\begin {lemma}
  {     Let $ Y \in { \mathcal S ^ q } $ be irreducible with $ | Y | >
    1 $. Then the pair $ ( M_{ \large{ * } ^ s Y } , \large{ * } ^ s Y
    ) $ is permutation isomorphic to the pair $ ( M_Y \wr {\rm Sym} (
    s ) , Y ^ { \wr s } ) $.  }  {      Recall that the subgroup $ T $
    satisfies $ T \cong {\rm Sym} ( s ) $.  By Lemma
    \arabicc{lemma:StructureNProd}{} we have $ N_{ * ^ s Y } =
    N_Y\cdot T \cong N_Y \wr {\rm Sym} ( s ) $ and by Lemma
    \arabicc{lemma:StructureSProd}{} $ S_{ \large{ * } ^ s Y } =
    \large{ * } ^ s S_Y $.  Thus $ N_ { \large{ * } ^ s Y } / S_{
      \large{ * } ^ s Y } = N_Y \cdot T / ( S_Y * S_Y * \dots * S_Y )
    \cong ( N_Y / S_Y ) \cdot T \cong M_Y \cdot T \cong M_Y \wr {\rm
      Sym} ( s ) $.  So $ M_{ \large{ * } ^ s Y } \cong M_Y \wr {\rm
      Sym} ( s ) $ and this wreath product can be written out as $ (
    M_Y * M_Y * \dots * M_Y ) \cdot {\rm Sym} ( s ) $; keeping track
    of the isomorphisms we see that each factor $ M_Y $ acts on the
    corresponding factor $ Y $ of $ \large{ * } ^ s Y $ while the
    subgroup isomorphic to $ {\rm Sym} ( s ) $ still permutes these
    factors.  Therefore the identification $ M_{ \large{ * } ^ s Y }
    \cong M_Y \wr {\rm Sym} ( s ) $ yields the isomorphism $ \large{ *
    } ^ s Y \cong Y ^ { \wr s } $.  }
\end {lemma}
\begin {block}
  {     We have determined the structure of the permutation $ k M_{ *
      ^ s Y } $-module $ k * ^ s Y $ and of the permutation $ k M_{
      \Delta ^ s Y } $-module $ k \Delta ^ s Y $ for irreducible $ Y $. Now we turn to
    finding conditions under which $ * ^ s Y $ or $ \Delta ^ s Y $ is
    a closed element of $ {\mathcal S } ^ q $. The following
    elementary lemma from group theory will be useful.  }  {{\bf
      Irreducibles and Control Of Closure} }
\end {block}
\mylabel{lemma:GroupTheory}
\begin {lemma}
  {     Let $ \alpha $ be a finite set and let the subsets $ \alpha_i
    $ partition $ \alpha $. For each $ i $ let $ H_i $ be a subgroup
    of $ {\rm Sym} ( \alpha_i ) \leq {\rm Sym} ( \alpha ) $ that does
    not fix an element of $ \alpha_i $. Then $ C_{ {\rm Sym} ( \alpha
      ) } ( H_1 \times H_2 \times \dots \times H_t ) =  C_{ {\rm
        Sym} ( \alpha_1 ) } ( H_1 ) \times C_{ {\rm Sym} ( \alpha_ 2 )
    } ( H_2 ) \times \dots \times C_{ {\rm Sym} ( \alpha_ t ) } ( H_t
    ) .$ } {     Let $ g \in {\rm Sym} ( \alpha ) $ be any element that
    centralizes $ H_i $. Then the action of $ <g> $ on $ \alpha $ will
    admit as an invariant subset the set of elements of $ \alpha $ not fixed
    by $ H_i $; by assumption this set is
    $ \supp \ H_i = \alpha_i $, so
    $ g \in \Sym (\alpha_i) \times \Sym ( \alpha \setminus \alpha_i) $.
    This holds for each $ i $ and so
    $ g \in \Sym (\alpha_1) \times \Sym (\alpha_2) \times \dots \times \Sym (\alpha_t) $,
    from which it follows
    $ C_{ {\rm Sym} ( \alpha
      ) } ( H_1 \times H_2 \times \dots \times H_t ) \leq C_{ {\rm
        Sym} ( \alpha_1 ) } ( H_1 ) \times C_{ {\rm Sym} ( \alpha_ 2 )
    } ( H_2 ) \times \dots \times C_{ {\rm Sym} ( \alpha_ t ) } ( H_t
    ) .$ The reverse inclusion is trivial.
  }
\end {lemma}
\vspace*{0.1cm}
\par If $ X \in {\mathcal S } ^ q $ we define the {\it closure}, written $ c ( X ) $ 
as the set $ \Fix (Q_X ) $ for a choice $ Q_X $ of Sylow $ p $-subgroup
of $ \Stab(X) $. Closure is well defined up to equivalence.
\mylabel{corollary:SubMultiplicative}
\begin {corollary}
  {     {\it Submultiplicativity of Closure.} Assume $ X_1 , X_2 , \dots
    , X_t $ belong to $ {\mathcal S } ^ q$, and are exact. Then $    c (
    X_1 * X_2 * \dots * X_t ) \leq c ( X_1 ) * c ( X_2 ) * \dots * c (
    X_t ) $.  }  { Set $ X := X_1 * X_2 * \dots * X_t $,
    and let $ Q_i := Q_{ X_i } $,
    and define also
    $ Q = Q_1 \times \dots Q_t $.
    Then the $ p $-group $ Q $ is contained in $ S_X $, so we can
    take $ Q \leq Q_X $.
    We assume the $X_i $ are exact, so we can apply
    \arabicc{lemma:GroupTheory}{} with
    $ H_i = Q_i $, and get, taking $ G = G_X $ and $ G_i = G_{ X_i } $,
    that
    $$ H := C_{ G } ( Q ) = C_{ G_1 } ( Q_1 ) \times \dots C_{ G_t } ( Q_t ) .$$
    Let $ \Xi_i = \Xi_{ X_i } $;
    then we have $ \Xi_1 * \dots * \Xi_t \subseteq \Xi_X $, and therefore
    $$ H \cap (  \Xi_1 * \dots * \Xi_t ) \subseteq H \cap \Xi_X .$$
    We claim that equality holds. Namely taking $ x \in \Xi_X $
    that commutes with $  Q = Q_1 \times \dots Q_t $, then each orbit
    of $ x $ is contained in $ X_i $ for some $ i $ and therefore
    $ x \in \Xi_1 * \dots * \Xi_t $. So $ H \cap \Xi_X $ is contained
    in $ H \cap (  \Xi_1 * \dots * \Xi_t ) $. This last set is equal to
    $ ( C_{ G_1 } ( Q_1 ) \cap \Xi_1 ) * \dots * ( C_{ G_t } ( Q_t ) \cap \Xi_t ) $.
    Hence
    $$ \Fix_{ \Xi_1 } ( Q_1 ) * \dots * \Fix_{ \Xi_t } ( Q_t ) 
       = C_G ( Q ) \cap \Xi_X \supseteq C_G ( Q_X ) \cap \Xi_X = {\rm Fix}_{ \Xi_X } ( Q_X  ), $$
    which is what we needed.
    }
\end {corollary}
\mylabel{corollary:ProductIsClosed}
\begin {corollary}
  {     Let $ X_1 , X_2 , \dots , X_t $ be exact closed members of $
    {\mathcal S } ^ q$. Then $ X_1 * X_2 * \dots * X_t $ is exact and closed.  }
  {  
    Write $ X = X_1 * X_2 * \dots * X_t $.
    We have $ c ( X_1 * X_2 * \dots * X_t ) \leq c ( X_1 ) * c (
    X_2 ) * \dots * c ( X_t ) = X_1 * X_2 * \dots * X_t $; thus we
    must have $ c ( X_1 * X_2 * \dots * X_t ) = X_1 * X_2 * \dots *
    X_t $ and $ X $ is closed.
    \par
    Let $ \alpha_j = \supp X_j $, writing 
    and $ \alpha = \bigcup \alpha_j $ and $ Q_j = Q_{ X_j } $.
    Then by assumption $ Q_j $ is fixed point free on $ \alpha_j $.
    We are free to assume $ Q_1 * Q_2 * \dots * Q_t \leq Q_X $
    and so $ Q_X $ is fixed point free on $ \alpha $, so $ X $ is
    exact.
}
\end {corollary}
\mylabel{lemma:DeltaProductIsClosed}
\begin {lemma}
  {     Let $ X $ be a closed exact element of $ {\mathcal S } ^ q $.
    If $ i $ is a power of $ p $ then $ \Delta ^ i X $ is
    closed.  \newline 
    When $ | X | > 1 $, then $  \Delta ^ i X $ is closed if and only if
    $ i $ is a power of $ p $.
    }  
    { By \arabicc{corollary:ProductIsClosed}{} $ \large{ * } ^ i X $ is
    closed for all $ i > 1 $.
    Moreover by definition
    $ \Delta ^ i X \subset \large{ * } ^ i X $.
    Therefore we always have
    $
    c ( \Delta ^ i X ) \subset c ( \large{ * } ^ i X ) = \large{ *
    } ^ i X .$
    In other words, letting $ \Xi = \Xi_{ \Delta ^ i X } $,
    we know that
    $$ \Fix_{ \Xi} ( Q_{ \Delta ^ i X } ) \subseteq ( \large{ * } ^ i X ). $$

    We must determine when the fixed point are only the diagonal.
    We take $ Q_{ \Delta ^ i X } = Q_X \wr Q $
    where $ Q $ is a Sylow $ p $-subgroup of $ \Sym ( i ) $.

    (a) Assume that $ i $ is a power of $ p $; then $ Q $
    is transitive on the support sets of the base group, and it
    follows that $ \Fix_{ * ^ i X } = \Delta ^ i X $ and therefore
    the fixed points of $ Q_{ \Delta ^ i X } $ in $ \large{ * } ^ i X $
    are as well just $ \Delta ^ i X $.

   (b) Assume $ | X | > 1 $. Suppose $ i $ is not
    a power of $ p $; then $ Q $ is not
    transitive on the support sets of the base group.
    In this case it has fixed points which which belong to
    $ \Delta ^ r X * \Delta ^ {i-r} X $ for some $ 1 \leq r \leq i $
    but not to $ \Delta ^ i X $. So $ \Delta ^ i X $ is not closed.}
\end {lemma}

    As a partial converse to the above we have the following.
\mylabel{projectivedeltanotclosed}
\begin{lemma}
{
    Let $ X $ be an irreducible closed projective element of $ {\mathcal S } ^ q $.
    If $ i \geq p $ is a power of $ p $ then $ \Delta ^ i X $ is
    not closed.     
}   
{
    By assumption $ X = \Xi_X $, and since the group
    $ \langle \Xi_X \rangle $ is transitive on $ \supp X $
    it follows that $ X $ is irreducible.
    \par Let $ S = S_{ \Delta ^ i X } $ and
    $ Q = Q_{ \Delta ^ i X } $,
    and write $ X_\Delta = \Delta ^ i X = \Delta X_1 \dots \Delta X_i $
    where $ \supp X_j = \alpha_j $ for disjoint $ \alpha_j $
    and $ \alpha = \bigcup \alpha_j $.
    We can now apply \arabicc{lemma:StructureS}{}
    and get
    $$ S =  ( S_{X_1} * \dots * S_{X_i} ) \cdot T,$$
    where $ T \cong {\rm Sym} ( i ) $. By \arabicc{section:FurtherWrProd}{}
    $ T $ is not transitive on $ \alpha $.
    Taking Sylow $ p $-sugroups we get
    $$ Q = ( Q_{ X_1 } * \dots * Q_{ X_i } ) \cdot P = P ,$$
    for a Sylow $p $-subgroup $ P $ of $ T $, since
    $ Q_X = 1 $ by assumption.
    Then $ X_\Delta $ is irreducible by
    \arabicc{lemma:DeltaIsIrreducible}{}
    but not exact, since $ P \leq T $ is not transitive;
    nor is $ X_\Delta $ projective; so $ X_\Delta $ cannot be closed by
    \arabicc{lemma:ExactOrProjective}{}.  
}
\end{lemma}
\begin {remark}
  {     We have gathered enough information to give sufficient
    conditions that given an irreducible fixed point set $ Y \in
    {\mathcal S } ^ q $ the element $ * ^ i Y $ or the element $
    \Delta ^ i Y $ is also a fixed point set. These conditions are set
    out in the next theorem.  }  { }
\end {remark}
\markright{Coprime Sets}
\mylabel{lemma:CoprimeSets}
\begin {theorm}
  {     Recall the definition of $ {\mathcal S} ^ q $ in
    \arabicc{definition:FS}{}.  Let $ Y \in { \mathcal S} ^ q $ be an
    exact (in the sense of \arabicc{block:Exact}{}) irreducible fixed
    point set. Let $ s < p $ be an integer and let $ t $ be a power of
    $ p $.  Then
        \begin {enumerate}
        \item $ * ^ s Y $ is a fixed point set. \\
        \item $ \Delta ^ t Y $ is an irreducible fixed point set. \\
        \end {enumerate}        
      } { Recall that $ Y $ is closed, as it is a fixed point set.
        By \arabicc{corollary:ProductIsClosed}{}
        the element $ * ^ s Y $ is closed. We want to show that
        in addition the permutation $ k M_ { * ^ s Y } $-module $ k *
        ^ s Y $ admits a projective summand. \par
        Assume first that $ | Y | > 1 $. Then by 
        \arabicc{lemma:StructureOfProduct}{}  the pair $ ( M_{
          \large{ * } ^ s Y } , \large{ * } ^ s Y ) $ is
        permutation isomorphic to the pair
        $ ( M_Y \wr {\rm Sym} ( s ) , Y ^ { \wr s } ) $. Thus the
        existence of the projective summand is a consequence of
        Lemma \arabicc{lemma:WrProdHasProj}{}.
        \par
        Now assume that $ | Y | = 1 $. Then
        $  \large{ * } ^ s Y = \Delta ^ s Y $
        and by \arabicc{lemma:StructureOfDeltaProduct}{},
        $ ( M_{ \Delta ^ s Y } , \Delta ^ s Y ) $ is permutation
        isomorphic to $ ( M_Y, Y ) $. Hence by assumption the
        permutation module $ k \Delta ^ s Y $ is projective as a module for
        $ M_{ \Delta ^ s Y } $. So $ \large{ * } ^ s Y $
        is also a fixed point set in this case.
        \par 
        That $ \Delta ^ t Y $ is irreducible comes from Lemma
        \arabicc{lemma:DeltaIsIrreducible}{}.   By Lemma
        \arabicc{lemma:DeltaProductIsClosed}{} the element $ \Delta ^
        t Y $ is closed, and we want to show that       the
        permutation $ k M_ { \Delta ^ t Y } $-module $ k \Delta ^ t Y
        $ admits a projective   summand. By Lemma
        \arabicc{lemma:StructureOfDeltaProduct}{}       the     pair $
        ( M_{ \Delta ^ t Y } , \Delta ^ t Y ) $ is permutation
        isomorphic to the pair $ ( M_Y , Y ) $. The existence of the
        projective summand follows by assumption{}.  }
\end {theorm}
\setcounter{chapterno}{5}
\begin {nblock}
  {     We have examined the structure afforded by products such as $
    * ^ s Y $, where $ Y \in {\mathcal S } ^ q $.  Let $ X , Y \in
    {\mathcal S } ^ q $ with possibly $ X \neq Y $. In this section we
    consider the product $ X * Y $ and the permutation $ k M_{ X * Y }
    $-module it affords. We also introduce the notion of `coprime'.
    
  } {Structure of Products of Distinct Sets.}
\end {nblock}
\mylabel{properties:mxy}
\begin {block}
  {     In order to determine the permutation structure of $ X * Y $
    we may assume as usual that $ X $ and $ Y $ have disjoint support
    and that $ X * Y = X \times Y $ . Then the subgroups $ N_X * N_Y
    $, $ S_{ X * Y }$ and $ S_X * S_Y $ are all contained in $ N_{ X *
      Y } $.  Thus we have the diagram of subgroups
                
        \begin{picture}(10,0)(-60,10) 
        \put( 57 , 3){$ N_{ X * Y } $ }
        \put (30,-20){ \line(1, 1){20}}   
        \put( 79 , 0){ \line(1,-1){20} }   
        \put ( 101 , -29){ $ N_X * N_Y $ }
        \put( 20 , -29) { $ S_{ X * Y }$ } 
        \put (30,-34){ \line(1, -1){20}}   
        \put( 99 , -34){ \line(-1,-1){20} }   
        \put(50 , -63){ $ S_X * S_Y  $. } 
        \end{picture} 
        \vspace{3cm}    \newline        We also have $ S_X * S_Y = S_{
          X * Y } \cap ( G_X * G_Y ) = S_{ X * Y } \cap ( N_X * N_Y )
        $.      The direct product $ M_X \times M_Y $ therefore
        satisfies       $ M_X \times M_Y = \dfrac{N_X}{S_X} \times
        \dfrac{N_Y}{S_Y} \cong \dfrac{N_X \times N_Y}{S_X \times S_Y}
        = \dfrac{N_X * N_Y}{S_X * S_Y} = \dfrac{ N_X * N_Y }{S_{ X * Y
          } \cap ( N_X * N_Y )} \cong   \dfrac{(N_X * N_Y) \cdot S_{X
            * Y} }{S_{X*Y}} \leq \dfrac{N_{X * Y }}{S_{X*Y}} =: M_{ X
          * Y}$.        This means that the direct product $ M_X
        \times M_Y $ can be naturally identified with   a subgroup of
        $ M_{ X * Y } $.  }      {{\bf Properties of $ M_{ X * Y } $}
      }
\end {block}
\begin {block}
  {     Let $ X $ and $ Y $ be as in the previous section. The groups
    $ M_X $ and $ M_Y $ are factor groups of subgroups of $ {\rm Sym}
    ( {\mathbb N} ) $ and as yet the expression $ M_X * M_Y $ is not
    defined.  Therefore we shall use $ M_X * M_Y $ to denote the
    subgroup of $ M_{ X * Y } $ which is naturally identifiable with
    the direct product $ M_X \times M_Y $. This subgroup is described
    in the previous section.  }  {{\bf Convention} }
\end{block}
\begin {block}
  {     Let $ X $ and $ Y $ be two elements of $ {\mathcal S } ^ q $.
    Then up to equivalence $ X \subset \Xi_{qm} ^ q $ and $ Y \subset
    \Xi_{qn} ^ q $ for some positive integers $ m $ and $ n $. Now
    suppose that $ X $ and $ Y $ have disjoint support, so that we may
    write $ X * Y = X \times Y $.  One way to view $ X * Y $ is that
    it is a subset of the $ {\rm Sym} ( q m ) * {\rm Sym} ( q n )
    $-space $ \Xi_{qm} ^ q * \Xi_{qn} ^ q $. The subgroup $ N_X * N_Y
    $ acts naturally on the Cartesian product $ X * Y $ and affords
    the permutation $ k ( N_X * N_Y ) $-module $ k X \otimes k Y $.  A
    second point of view is that $ X * Y $ is a subset of the $ {\rm
      Sym}_{ q ( m + n ) } $-space $ \Xi_{ q ( m + n ) } ^ q $. In
    this case $ X * Y $ is a permutation $ N_{ X * Y} $-space and
    affords the permutation $ k N_{ X * Y } $-module $ k ( X * Y ) $.
    There is at least as much structure in the second instance as in
    the first, and indeed we have an isomorphism of $ k N_X * k N_Y
    $-modules
    $$
    {\rm Res} ^ { N_{ X \times Y } } _ { N_{ X } * N_{ Y } } ( k (
    X * Y ) ) \cong k X \otimes k Y .$$
    Quotienting out the
    stabilizers of the various actions produces an isomorphism of $ k
    ( M_X * M_Y ) $-modules
    $$
    {\rm Res} ^ { M_{ X \times Y } } _ { M_{ X } * M_{ Y } } ( k (
    X * Y ) ) \cong k X \otimes k Y .$$
    
    The second point of view is the one we have to work with in order
    to show that $ X * Y $ is a fixed point set. However, we will be
    very interested in those cases where the second view adds no more
    structure to the first, which leads to the following definition.
    \newline {\bf Definition.} 
    We say that $ X $ and $ Y
    $ are { \bf coprime} if $ N_{ X * Y } = N_X * N_Y $, or
    equivalently if $ N_{ X * Y} \leq G_X * G_Y $.  }  {{\bf Coprime
      Sets} }
\end {block}
\mylabel{proposition:ProdBehaviour}
\begin {block}
  {      Let $ X $ and $ Y $ be coprime elements of $ {\mathcal S } ^
    q $. Then
        \begin {enumerate}
        \item we have $ S_{ X * Y } = S_X * S_Y $
        \item we have $ M_{ X * Y } = M_X * M_Y $.
        \end {enumerate}        
        {\bf Proof.} By hypothesis we have $ N_{ X * Y} \leq G_X * G_Y
        $ and moreover $ N_{ X * Y} = N_X * N_Y $.      \newline (1)
        Since $S_X * S_Y = S_{ X * Y } \cap ( G_X * G_Y ) $ it follows
        that $ S_{ X * Y } = S_X * S_Y $.       \newline (2) This
        holds because $ M_X * M_Y $ is a subgroup of $ M_{ X * Y} $
        and     $ M_{ X * Y} := N_{ X * Y} / S_{ X * Y} = ( N_X * N_Y
        ) / ( S_X * S_Y )       \cong ( N_X / S_X ) \times ( N_Y / S_Y
        ) \cong M_ X * M_Y $ (see \arabicc{properties:mxy}{}). \qed } {{\bf Proposition} }
\end {block}
\mylabel{block:StructCoprime}
\begin {block}
  {     Thus if $ X $ and $ Y $ are coprime we must have $ M_{ X * Y }
    = M_{ X } * M_{ Y } $. In this case the isomorphism of permutation
    $ k ( M_X * M_Y ) $-modules
    $$
    {\rm Res} ^ { M_{ X \times Y } } _ { M_{ X } * M_{ Y } } ( k (
    X * Y ) ) \cong k X \otimes k Y .$$
    is in fact an isomorphism of
    permutation $ k M_{ X \times Y } $-modules
    $$
    k ( X * Y ) \cong k X \otimes k Y .$$
  } {{\bf Structure of
      Products of Coprime Sets} }
\end {block}
\mylabel{lemma:ClosedProductIsFixed}
\begin {lemma}
  { Suppose $ Y $ and $ Z $ are two coprime elements of 
    $ {\mathcal S } ^ q $ which are fixed point sets.
    Assume $ Y * Z $ is closed.
    Then $ Y * Z $ is
    a fixed point set.  } 
{   Recall that fixed point sets are closed.
    By assumption $ Y * Z $ is closed, so it
    is sufficient to show that the permutation
    $ k M_{ Y * Z } $-module $ k ( Y * Z ) $ has a projective summand.
    The fact that $ Y $ and $ Z $ are coprime implies that
    $ M_{ Y * Z } = M_Y \times M_Z $
    and that $ k ( Y * Z ) $ is isomorphic to $ k
    Y \otimes k Z $. Thus it suffices to show that the
    $ k M_Y $-module $ k Y $ and the $ k M_Z $-module
    $ k Z $ each admit projective
    summmands. This is known since $ Y $ and $ Z $ are
    fixed point sets.}
\end {lemma}
\mylabel{lemma:CoprimeProductIsFixed}
\begin {lemma}
  { Suppose $ Y $ and $ Z $ are two coprime exact elements of $
    {\mathcal S } ^ q $ which are fixed point sets.  Then $ Y * Z $ is
    a fixed point set.  } 
{   
    By Corollary
    \arabicc{corollary:ProductIsClosed}{} we know that $ Y * Z $ is
    closed, and we may apply \arabicc{lemma:ClosedProductIsFixed}{}.}
\end {lemma}
\begin {block}
  {      {\it Let $ X , Y \in {\mathcal S } ^ q $. Then $ X $ and $ Y
      $ are coprime if and only if they do not have an irreducible
      factor in common.} \newline {\bf Proof.}  We may assume that $ X
    $ and $ Y $ have disjoint supports so that $ X * Y = X \times Y $.
    \par Suppose first that $ X $ and $ Y $ share the irreducible
    factor $ Z $. By hypothesis we have $ X = Z_X \times W_X $ and $ Y
    = Z_Y \times W_Y $ where $ Z_X $ and $ Z_Y $ are equivalent to $ Z $.
    Let $ \alpha = \supp Z_X $ and $ \beta = \supp Z_Y $.

    Recalling the definition made in Section
    \arabicc{section:FurtherWrProd}{} let $ t := g_{ \alpha \supp \ 
      Z_Y } \in {\rm Sym} ( {\mathbb N} ) $ be the unique
    order-preserving involution that maps $ \alpha $ onto $ \supp \ 
    Z_Y $ and satisfies $ \supp \ t = \alpha \cup \beta $.
    \par We may assume that $ Z_Y = Z_X ^ t $. Now $ X ^ t = Z_Y
    \times W_X $ and $ Y ^ t = Z_X \times W_Y $ so $ ( X \times Y ) ^
    t = X ^ t \times Y ^ t = X \times Y $. It follows that $ t \in N_{
      X * Y } $. However, $ \alpha \subset \supp \ X $ while
    $ \alpha  ^ t = \beta \subset \supp \ Y $ so it is clear
    that $ t $ does not leave $ \supp \ X $ or $ \supp \ Y $ invariant.
    Thus $ t \notin G_X \times G_Y $, and $ X $ and $ Y $ are not coprime.
        \par 
        Suppose instead that $ X $ and $ Y $ are not coprime and let $
        t \in N_{ X * Y } $ with $ t \notin G_X * G_Y $.        Let $
        X = X_1 \times X_2 \times \dots \times X_s $ and $ Y = Y_1
        \times Y_2 \times \dots \times Y_r $ be decompositions  into
        irreducibles. Then $ X \times Y = X_1 \times X_2 \times \dots
        \times X_s \times Y_1 \times Y_2 \times \dots \times Y_r $
        is a decomposition into irreducibles. If for every $ 1 \leq i
        \leq s $ there is $ j $ such that $ X_i ^ { \ t} = X_j $ then
        $ t $ leaves $ \supp \ X $ and hence $ \supp \ Y $ invariant. By
        assumption this is not the case. Therefore by   Corollary
        \arabicc{corollary:PermutesFactors}{} for some $ 1 \leq i \leq
        s $ and $ 1 \leq j \leq r $ we have $ X_i ^ { \ t } = Y_j $.
        Thus $ X_i = Y_j $ modulo equivalence and $ X $ and $ Y $
        share the irreducible factor $ X_i. \qed$
        
      } {{\bf Proposition} }
\end {block}
\mylabel{theorem:ProdIsFix}
\begin {theorm}
  {     Suppose that $ Y_1 , Y_2 , \dots , Y_t \in {\mathcal S } ^ q $
    is a list of pairwise distinct irreducible exact fixed point sets
    and let $ s_1 , s_2 , \dots , s_t $ be positive integers
    satisfying $ s_i < p $ for every $ i $. Then $ Y_1 ^ { s_1} * Y_2
    ^ { s_2 } * \dots * Y_t ^ { s_t } $ is a fixed point set.  }
  {     This follows directly from Theorem \arabicc{lemma:CoprimeSets}{}
    and Lemma \arabicc{lemma:CoprimeProductIsFixed}{} because distinct
    irreducible elements of $ {\mathcal F} $ are coprime.}
\end {theorm}
\vspace*{2mm}
\par To show how the assumption of exactness is necessary in the
previous theorem we have the following.
\mylabel{lemma:projectiveproductnotclosed}
\begin{lemma}
{
    Let $ Y $ and $ Z $ be closed projective elements
    of $ {\mathcal S } ^ q $. Then $ Y * Z $ is not closed.
}
{
    Write $ X = Y * Z $.   
    By assumption we have $ Y = \Xi_Y $ and $ Z = \Xi_Z $.
    In
    particular $ Y $ and $ Z $ are irreducible. If $ Y \neq Z $
    then $ Y $ and $ Z $ are coprime and we have
    $ S_X = S_Y * S_Z $ by
    \arabicc{proposition:ProdBehaviour}{} . Otherwise $ Y = Z $
    and again $ S_X = S_Y * S_Z $ by
    \arabicc{lemma:StructureSProd}{}. Accordingly we may
    always choose
    $ Q_X = Q_Y * Q_Z $ : but by assumption this gives $ Q_X = 1 $.
    So if $ X $ is closed then $ X = \Xi_X $ which is a contradiction
    since $ \Xi_X $ is irreducible.}
\end{lemma}


\setcounter{chapterno}{6}
\setcounter{lastplace}{1}

\markboth{Reducing Fixed Point Sets}{Reducing Fixed Point Sets}
\begin {nblock}
  {     We have seen how to compose fixed point sets from irreducible
    ones, and in this section we will confirm that such compositions
    are the general form of fixed point sets. The results in this
    section will tend to be the converses of those in the previous
    section. By the end of the section the problem will have been
    reduced to finding the irreducible fixed point sets, which we will
    do in the next and final section.  }  {Decomposing The Fixed Point
    Sets.  }
\end {nblock}
\begin {block}
  {     Given a fixed point set $ Z $ consider a decomposition of $ Z
    $ into irreducible factors - our goal is to show that each of the
    factors is again a fixed point set. This reduces the problem to
    determining all the irreducible fixed point sets.  }  {{\bf Divisors of Fixed Point Sets }}
\end {block}
\begin {block}
  {     If $ X \in {\mathcal S } ^ q $ let $ \tilde{Q} _ X := < Q_X ,
    X > $ be the subgroup of $ G_X $ generated by $ Q_X $ and $ X $.
    See Section \arabicc{notation:QSMagain}{} for the definitions of $
    Q_X $ and $ X $ and Section \arabicc{definition:FS}{} for the
    definition of $ {\mathcal S } ^ q $.  }  {{\bf Notation} }
\end {block}
\begin {block}
  {     In Section 3 we used Lemma \arabicc{lemma:FixCombo}{} to split
    the task of showing that a set $ Y $ is a fixed point set into two
    parts. Thus to show that $ Y $ is a fixed point set we show that
    (i) $ Y $ is closed, and that (ii) the permutation $ k M_Y
    $-module $ Y $ admits a projective summand. Since Lemma
    \arabicc{corollary:ProdHasProj}{} implies that this second
    property passes readily to any factors in a decomposition of $ Y
    $, the first results of this chapter deal with whether or not such
    factors are closed.  }  {{\bf Control Of Closure} }
\end {block}
\begin {remark}
  {     Suppose that $ Y $ and $ Z $ are elements of $ { \mathcal S }
    $ such that $ Y * Z $ is closed. We shall see in the proof of the
    next lemma that if $ S_{ Y * Z } = S_Y * S_Z $ then we may infer
    that $ Y $ and $ Z $ are closed. A sufficient condition for $ S_{
      Y * Z } = S_Y * S_Z $ to hold is that the subgroups $ < Y * Z >
    $ and $ < Y > * < Z > $ are equal. However, we shall take another
    approach and split this problem into two cases.  }  { }
\end {remark}
\mylabel{lemma:FactorsAreClosed}
\begin {lemma}
  {     Suppose that $ Y $ and $ Z $ are coprime elements of $ {
      \mathcal S } $ such that $ Y * Z $ is closed. Then $ Y $ and $ Z
    $ are closed.  }  { The coprime condition gives $ Q_{ Y * Z } =
    Q_Y * Q_Z $.  Thus $ Y * Z = c ( Y * Z ) = \Fix ( Q_{ Y * Z } ) =
    \Fix ( Q_Y * Q_Z ) \supset \Fix ( Q_Y ) * \Fix ( Q_Z ) $, giving $
    \Fix ( Q_Y ) = Y $ and $ \Fix ( Q_Z ) = Z $.  }
\end {lemma}
\markboth{The Irreducible Fixed Point Sets}{The Irreducible Fixed
  Point Sets}
\mylabel{lemma:SameFactorsAreClosed}
\begin {lemma}
  {     Suppose $ Y $ is an irreducible element of $ {\mathcal S } ^ q
    $ such that $ Y ^ i $ is closed.  Then $ Y $ is closed.  }
  {     If $ | Y | = 1 $ then $ Y ^ i = \Delta ^ i Y $. This case is
    dealt with in Lemma (7.10): no circularity of argument will result
    from this expedience.  By Lemma \arabicc{lemma:StructureSProd}{}
    we have $ S_{ Y ^ i } = S_Y ^ i $, and so $ Q_{ Y ^ i } = Q_Y *
    Q_Y * \dots * Q_Y $.  The proof is now the same as the proof of
    the last lemma.  }
\end {lemma}
\mylabel{corollary:FactorsClosed}
\begin {corollary}
  {     Suppose that $ Y $ and $ Z $ are coprime elements of $ {
      \mathcal S } $ such that $ Y * Z $ is a fixed point set. Then $
    Y $ and $ Z $ are fixed point sets.  }  {   The coprime condition
    gives $ M_{ Y * Z } = M_{ Y } \times M_{ Z } $ and $ k ( Y * Z )
    \cong k Y \otimes k Z $, by Section
    \arabicc{block:StructCoprime}{}.  By hypothesis the permutation $
    k M_{ Y * Z } $-module $ k ( Y * Z ) $ admits a projective
    summand. Thus by Corollary \arabicc{corollary:ProdHasProj}{} the
    permutation $ k M_Y $-module $ k Y $ and the permutation $ k M_Z
    $-module $ k Z $ each admit projective summands. By Lemma
    \arabicc{lemma:FactorsAreClosed}{} $ Y $ and $ Z $ are closed so
    Lemma \arabicc{lemma:FixCombo}{} implies that $ Y $ and $ Z $ are
    fixed point sets.  }
\end {corollary}
\par 
Note that by Lemma \arabicc{lemma:UniqueDecomp}{} the decomposition in
the statement of the following theorem is unique.
\mylabel{theorm:ReductionToIrreducibles}
\begin {theorm}
  {     Let $ Y \in {\mathcal S } ^ q $ be a fixed point set and write
    $ Y = Y_1 ^ { a_1 } * Y_2 ^ { a_2 } * \dots * Y_t ^ { a_t } $
    where the $ Y_i $ are pairwise distinct irreducible elements and
    the $ a_ i $ are positive integers. Then each $ Y_i $ is a fixed
    point set.  }  {    The list $ Y_1 ^ { a_1 } , Y_2 ^ { a_2 } ,
    \dots , Y_t ^ { a_t } $ is a list of pairwise coprime elements of
    $ {\mathcal S } $ so Theorem \arabicc{theorem:ProdIsFix}{} implies
    that each $ Y_i ^ { a_i } $ is a fixed point set. In particular $
    Y _i ^ { a_i } $ is closed for every $ i $.
    \newline Assume first $ | Y_i | > 1 $ for a given $ i $.  Now Lemma
    \arabicc{lemma:StructureOfProduct}{} tells us that for each $ i $
    $ k Y_i ^ { a_i } \cong Y_i ^ { \wr a_i } $ under the
    identification $ k M_{ Y_i ^ { a_i } } = k M_{ Y_i } \wr {\rm Sym}
    ( a_ i ) $. Hence Lemma
    \arabicc{lemma:DirectProductProjectiveSummand}{} asserts that the
    permutation $ k M_{ Y_i } $-module $ Y_i $ admits a projective
    summand. 
    \newline Now assume $ | Y_i | = 1 $;
    then $ Y_i^{a_i}= \Delta ^ { a_i } Y_i $ 
    and we can apply \arabicc{lemma:StructureOfDeltaProduct}{} which
    shows that the relevant permutation module has a projective summand.
    So in both cases, by \arabicc{lemma:FixCombo}{} each $ Y_i $ is a fixed point set.

Meanwhile Lemma \arabicc{lemma:FactorsAreClosed}{} and
    Lemma \arabicc{lemma:SameFactorsAreClosed}{} tell us that each $
    Y_i $ is closed, so by Lemma \arabicc{lemma:FixCombo}{} $ Y_i $ is
    a fixed point set.  }
\end {theorm}

\setcounter{chapterno}{7}
\begin {nblock}
  {     In this section we shall determine the structure of the
    irreducible fixed point sets (in the context of the present
    problem) and in doing so we will finally complete our description
    of a general fixed point set. The results of the previous sections
    will combine with the following propositions in order to arrive at
    the main theorem.  }  {Determining the Irreducible Fixed Point
    Sets.}
\end {nblock}
\begin {block}
  {     We will show in Corollary (7.6) that if $ Y $ is an
    irreducible fixed point set then $ \tilde{Q} $ must be
    transitive on $ \supp \ Y $. Thus $ Q_Y := <Q_Y, Y> $ transitively permutes the
    orbits on $ \supp \ Y $ of $ < Y > $. It will follow that $ Y =
    \Delta ^ i X $ for some $ i $ and $ X \in {\mathcal S } ^ q $,
    where the size of $ \supp \ X $ is equal to the size of a $ < Y >
    $-orbit on $ \supp \ Y $. In Theorem (7.14) we find a bound for the
    orbit size, reducing the problem to the study of a finite, and in
    fact a very small, number of cases.  }  {{\bf Reducing the
      Irreducible Fixed Point Sets} }
\end {block}
\mylabel{lemma::transitiverootclosed}
\begin {lemma}
  {     Suppose that $ Y \in {\mathcal S } ^ q $ is closed irreducible
    Suppose further that $ Y = \Delta ^ s Z $ for some $ Z \in
    {\mathcal S } ^ q $. Then $ Z $ is closed irreducible.  }  {  
    The congruence
    $ ( X_1 * X_2 ) \Delta ( X_1 * X_2 ) \cong ( X_1 \Delta X_1 ) * ( X_2 \Delta X_2 ) $
    holds generally for $ X_1, X_2 $ in $ {\mathcal S } ^ q $ so
    $ Z $ must be irreducible.
    \par By Lemma
    \arabicc{lemma:StructureS}{} we have $ S_Y = S_{\Delta ^ s Z } =
    S_Z \cdot T = * ^ s S_Z \cdot T \cong S_Z \wr {\rm Sym} ( s ) $
    where $ T \cong {\rm Sym} ( s ) $ is the subgroup defined in
    \arabicc{section:FurtherWrProd}{} that permutes the factors of the
    base subgroup $ * ^ s S_Z \leq S_Y $.  It follows from this that $
    Q_Y = * ^ s Q_Z \cdot R $ where $ R \leq T $ is a Sylow $ p
    $-subgroup of $ T $.  Let $ z \in c ( Z ) $ - then $ c \in \Fix \ 
    Q_Z $. Clearly $ y := z * z * \dots * z \in \Fix \ * ^ s Q_Z $;
    but $ y $ is also stabilized by $ T $.  Thus $ y \in \Fix \ * ^ s
    Q_Z \cdot T \subset \Fix \ Q_Y = Y $, and $ z \in Z $. So
    $ Z $ is closed.}
\end {lemma}

\begin {block}
  {     If $ G := {\rm Sym} ( m ) $ is the {\rm Sym}metric group on $
    m > 0 $ letters and $ H \leq G $ is a subgroup of $ G $ which
    leaves invariant the two subsets $ \alpha \subset \{ 1 , 2 , \dots
    , m \} $ and $ \beta \subset \{ 1 , 2 , \dots , m \} $, and if
    these subsets satisfy ( i) $ \alpha \cap \beta = \emptyset $ and
    (ii) $ \alpha \cup \beta = \{ 1 , 2 , \dots , m \} $, then we know
    that $ H \leq G $ sits inside the Young subgroup $ {\rm Sym} (
    \alpha ) * {\rm Sym} ( \beta ) \leq G = {\rm Sym} ( m ) $.
    However, it is not necessarily true that there is a subgroup $ H_
    \alpha \leq {\rm Sym} ( \alpha ) $ of $ {\rm Sym} ( \alpha ) \leq
    G $ and a subgroup $ H_ \beta \leq {\rm Sym} ( \beta ) $ of $ {\rm
      Sym} ( \beta ) \leq G $ for which $ H = H_ \alpha * H_ \beta $.
        \par 
        The next theorem finds conditions which do
        allow us make this assumption.  }  {{\bf Orbits and
          Factorizations} }
\end {block}

 \vspace*{2mm} { Let $ \tilde{Q}_X = \langle Q_X, X \rangle. $ }
\mylabel{theorem:VertexIsProduct}
\begin {theorm}
{   Assume $ X \in {\mathcal S } ^ q $ is exact and closed.  Assume
    the group $ \tilde{Q}_X $ has orbits
    $ \lambda_1 , \lambda_2 , \dots , \lambda_s $ 
    on the support of $ X $. 
    Then 
    \newline (1) $ X = X_1 * X_2 * \dots * X_s $ where $ X_i $ has
    support $ \lambda_i $ and is closed as a subset of $ \Omega_{ X_i } $
    and $ X_i $ exact.
    \newline (2) $Q_X = Q_1 * Q_2 * \dots * Q_s $ where $ Q_i = Q_{ X_i } $.
}
{    
     It suffices to consider $ s = 2 $. 
     Write $ Q $ and $ \tilde{Q} $ for $ Q_X $ and $ \tilde{Q}_X $. 
     We have $ \tilde{Q} \subseteq H_1 \times H_2 $
     where $ H_i := \Sym ( \lambda_i ) $. Let $ \pi_i $ be the projection
     onto $ H_i $, and set
     $$ Q_i := Q_{ \pi_i }, X_i := X_{ \pi_i }. $$
     Then we have
     \newline (i)  $ Q \subseteq Q_1 \times Q_2 $,
     \newline (ii) $ X \subseteq X_1 * X_2 $.
     \newline We claim that $ Q_i $ fixes each element of $ X_i $. Namely, take
     $ u_1 \in Q_1 $, then choose
     $ u_2 \in Q_2 $ such that $ u := u_1u_2 \in Q $.
     Next, take $ x_1 \in X_1 $, and take $ x_2 \in X_2 $ such that
     $ x := x_1 * x_2 $ belongs to $ X $. Then we have
      $$ x = x ^ u = ( x_1 ^ {u_1} ) * ( x_2 ^ {u_2} ) $$
      and $ x_1 = x ^ { u_1 } $. This shows that $ Q_1 \subseteq S_{X_1} $,
      and similarly $ Q_2 \subseteq S_{ X_2 } $.

     Now take $ R_i $ to be a Sylow $ p $-subgroup  of $ S_X $, containing
     $ Q_i $. Then $ R_1 \times R_2 $ fixes each element of $ X_1 * X_2 $
     and in particular it fixes each element of $ X $. So it follows
     that the order of $ R_1 \times R_2 $ is $ \leq $ the oreder of some Sylow
     $ p $-subgroup of $ S_X $ and hence $ \leq | Q | $.

     On the other hand, from $ Q \leq Q_1 \times Q_2 \leq R_1 \times R_2 $
     we have $ | Q | \leq | R_1 \times R_2 | $ and it follows that
     $$ Q= Q_1 \times Q_2 = R_1 \times R_2 $$
     and $ Q_i = R_i $ which we can take for $ Q_{ X_i } $. This completes the
     proof of part (2).
     
     It remains to show that $ X_1 * X_2 \subseteq X $.
     Let $ \Omega_i = \Xi_{ X_i } $, then
     $ \Omega_1 * \Omega_2 \subset \Xi_X $. We have
     $$ X_1 * X_2 \subseteq \Fix_{ \Omega_1 } ( Q_1 ) * \Fix_{ \Omega_2 } ( Q_2 )
        \subseteq \Fix_{ \Omega_1 * \Omega_2 } ( Q ) \subseteq 
        \Fix_{ \Xi_X } ( X );$$
     and $ {\rm Fix}_{ \Xi_X } ( Q ) = X $ since $ X $ is closed.
     So $ X = X_1 * X_2 $ and $ X_i = \Fix_{ \Omega_i } ( Q_i ) $,
     that is $ X_i $ is closed.     
}
\end {theorm}
\mylabel{corollary:QXPrimeIsTransitive}
\begin {corollary}
  {     Suppose that $ X \in {\mathcal S } ^ q $ is closed, exact and
    irreducible. Then $ < Q_X , X > $ is a transitive subgroup of $
    G_X = {\rm Sym} ( \supp \ X ) $.  }  {       This follows
    immediately from the previous Theorem
    \arabicc{theorem:VertexIsProduct}{}.  }
\end {corollary}
\mylabel{corollary:QXIsTransitive}
\begin {corollary}
  {     Suppose that $ X \in {\mathcal S } ^ q $ is closed, exact and
    irreducible. Let $ \lambda_1 , \lambda_2 , \dots , \lambda_t $ be
    the orbits on $ \supp \ X $ of $ < X > $. Then $ Q_X $ transitively
    permutes the $ \lambda_i $. In particular these orbits have the
    same size.  }  {    Firstly the subgroup $ < Q_X , X > \leq G_X $
    normalizes its subgroup $ < X > $ so $ < Q_X , X > $ permutes the
    $ \lambda_i $, and it must do so transitively since $ < Q_X , X >
    $ is transitive on $ \supp \ X $ (by Corollary
    \arabicc{corollary:QXPrimeIsTransitive}{}).  Now $ < X > $ acts
    trivially on the $ \lambda_i $ and it follows that $ Q_X $ must
    transitively permute the $ \lambda_i $.  }
\end {corollary}
\markright{Transitive Sets}
\begin {block}
  {     Let $ X \in {\mathcal S } ^ q $. We say that $ X $ is a {\bf
      transitive} set if $ < X > $ is transitive in its action on $
    \supp \ X $.  }  {{\bf Transitive Sets} }
\end {block}
\mylabel{corollary:FormOfIrreducibles}
\begin {corollary}
  {     Suppose that $ X \in {\mathcal S } ^ q $ is closed, exact and
    irreducible.  Then $ X = \Delta ^{p ^ i} Y $ for some
    integer $ i \geq 0 $ and some transitive set $ Y \in {\mathcal S } ^ q $.
  }  {  Let $ \lambda_1 , \lambda_2 , \dots , \lambda_t $ be the
    orbits of $ < X > $ on $ \supp \ X $ - then the $ \lambda_i $
    have equal sizes, by \arabicc{corollary:QXIsTransitive}{}.
    Put $ \lambda := \lambda_1 $.  Furthermore we have
    $$
    X \leq X_{ \lambda_1 } * X_{ \lambda_2 } * \dots * X_{
      \lambda_t } = * ^ t X_{ \lambda } .$$
    Now $ Q_X $ stabilizes $ X
    $ {\it and} transitively permutes the $ \lambda_i $. It follows
    that $ X \leq X_{ \lambda_1 } \Delta X_{ \lambda_2 } \Delta \dots
    \Delta X_{ \lambda_t } = \Delta ^ t X_{ \lambda } $, from which it
    follows that $ X = \Delta ^ t Y $ where $ Y := ( X ) \pi_{ \lambda
    } $. Also, we see that $ t $ must be a $ p $-power, for $ t $ is
    the size of the set $ \{ \lambda_1 , \lambda_2, \dots , \lambda_t
    \} $, and the elements of this set are transitively permuted by $
    Q_X $.  Now $ \lambda $ is an $ < X > $-orbit on $ \supp \ X $ so $
    < ( X ) \pi_{ \lambda } > $ acts transitively on $ \lambda $,
    whence $ Y $ is a transitive set.}
\end {corollary}
\mylabel{lemma:FixFreeElement}
\begin {lemma}
  {     Let $ X \in {\mathcal S } ^ q $ be closed and exact.  Then $ Z
    Q_X $ contains an element of order $ p $ fixed point free on $
    \supp \ X $.  }  {   By Lemma \arabicc{theorem:VertexIsProduct}{}
    it suffices to consider the case that $ \tilde{Q}_X := < Q_X , X
    > \leq G_n = {\rm Sym} ( 2 n ) $ is transitive in its action on $
    \supp \ X $. If $ u \in Z Q_X $ is an element of $ Z Q_X $ then $ u
    $ centralizes $ Q_X $ and, since $ u \in Q_X $, $ u $ centralizes
    $ X = \Fix ( Q_X ) $. Thus $ Q_X $ and $ X $ permute the fixed
    points on $ \supp \ X $ of $ < u > $ and in fact the transitive
    subgroup $ \tilde{Q}_X $ permutes these fixed points.  It follows
    that if we choose $ u \in Z Q_X $ so that $ u $ has order $ p $
    then $ u $ cannot fix any points in $ \supp \ X $.  }
\end {lemma}
\mylabel{corollary:StabP}
\begin {block}
  {   
      Now let $ Y $ be an irreducible exact fixed point set
      in $ {\mathcal S } ^ q $
      By \arabicc{lemma:FixFreeElement}{} there is some $ z $
      central in $ Q_Y $ which is fixed point free and has order
      $ p $.
      \par
      We take a Sylow $ p $-subgroup $ P $ of $ C_{ G_Y } ( z ) $
      which contains $ Q_Y $. Then $ P $ is even a Sylow $ p $-subgroup
      of $ G_Y $ ( as $ P $ is of the form $ C_p \wr R $ for a Sylow
      subgroup of $ \Sym ( s ) $ ).
      Then by \arabicc{lemma:vertexintersection}{} there is some $ y \in Y $
      such that $ Q_Y = \Stab_P ( y ) $. We will use these elements $ z $ 
      and $ y $ below.      
}
{
}
\end {block}
\markright{Determining The Transitive Fixed Point Sets}
\mylabel{block:StandardResult} \newline
\begin {block}
{      {\it Let $ Y $ be a closed element of $ {\mathcal S } ^ q $.
      Suppose that $ Q_Y $ contains subgroups $ R_1 , R_2 , \dots ,
      R_l $. Let $ \supp R_i = \alpha_i $
      and suppose that the $ \alpha_i $ are disjoint.
      Suppose also that $ * R_i \leq Q_Y $. Then $ Y \subseteq * {\rm Sym}
      ( \alpha_i ) $.}    
{     \newline {\bf Proof.}
      Apply
      \arabicc{lemma:GroupTheory}{} with $  H_i = R_i.\qed$}
}
{{\bf Proposition}}
\end {block}
\mylabel{convention::zconvention}
\begin {convention}
 {
    Let $ Y \in {\mathcal S } ^ q $ be an exact fixed point set.
    By Lemma \arabicc{lemma:FixFreeElement}{} there is an element $ z
    \in Z Q_Y $ which has order $ p $ and which is fixed point free on
    $ \supp \ Y $.  Taking the cycle decomposition of $ z $ we have $ z
    = z_1 * z_2 * \dots * z_s $ where each $ z_i $ is a $ p $-cycle.
    Let $ P_z $ be a Sylow $ p $-subgroup of $ C_{ G_Y } ( z ) $. Now
    $ C_{ G_Y } ( z ) $ certainly centralizes $ z $ and must permute
    the $ \{ z_i \} $.  Thus $ B_z := < z_1 > * < z_2 > * \dots * <
    z_s > $ is a normal $ p $-subgroup of $ C_{ G_Y } ( z ) $ - in
    fact $ B_z $ is the base subgroup of the wreath product $ C_{ G_Y
    } ( z ) \cong C_p \wr {\rm Sym} ( s ) $, as we have seen in
    \arabicc{block:ActOnProd}{} It follows that $ B_z \leq P_z $.
        \par By \arabicc{block:ActOnProd}{} we know that $ P_z $ is a Sylow $ p $-subgroup of $ G_Y $.
        We observe that in      in this context the roles of    $ p $
        and $ q $ have interchanged.  }  { }
\end {convention}
\vspace*{1mm}
\par The next theorem determines the possible transitive irreducible exact
fixed point sets. Recall that $ {\rm d} ( Y ) = | \supp \ ( Y ) | $.
\mylabel{lemma:BoundOnIrreducibles}
\begin {theorm}
{ 
    Assume that $ q $ is prime.
    Let $ Y \in {\mathcal S } ^ q $ be an irreducible fixed
    point set which is exact and transitive.
    Let $ {\rm d} ( Y ) = | \supp ( Y ) | $. Then one of the following
    holds
    \newline (1) $ {\rm d} ( Y ) = pq $
    \newline (2) $ {\rm d} ( Y ) = q  $.
    \newline In case (2), $ Y $ is the set of $ q $-cycles in 
    $ G_Y \cong \Sym ( q ) $.
}
{
    (a) Let $ z \in Z ( Q_Y ) $ be an element of order $ p $
    which is fixed point free on $ \supp ( Y ) $.
    [This exists by \arabicc{corollary:StabP}{}.]
    Then $ Q_Y \subseteq C_G ( z ) $. Take a Sylow $ p $-subgroup
    $ P := P_Z $ of $ C_G ( z ) $ containing $ Q_Y $.
    Then $ P \in Syl_p ( G ) $.
    [To see this: we know $ C_G ( z ) \cong C_p \wr \Sym ( s ) $
    where $ z $ has $ s $ $ p $-cycles, and $ z $ has support
    equal to the total support of $ G_Y $. So $ P $is of the form
    $ C_P \wr R $ where $ R $ is a Sylow subgroup of $ \Sym ( s ) $.
    On the other hand, a Sylow subgroup of $ G_Y $
    is exactly of this form, since $ p $ divides the size of its
    support.
    \par
    Now we use \arabicc{lemma:vertexintersection}{},
    which gives us that there is some $ y \in Y $ such that
    $ Q_Y = \Stab_P ( y ) $. Then $ z $ and $ y $ commute.
    \par
    (b) Since $ y $ and $ z $ commute, conjugation by
    $ y $ permutes the cycles of $ z $. Let $ I $ be an
    index set labelling the orbits, and write $ \overline{z}_i $
    for the product of the cycles in orbit $ i $.
    Then $ \overline{z}_i $ commutes with $ y $, for each $ i $.
    Recall from \arabicc{convention::zconvention}{} that $ P $
    contains each cycle of $ z $,
    so $ \overline{z}_i $ belongs to $ P $
    and therefore it belongs to $ \Stab_P ( y ) = Q_Y $.
    Now we apply \arabicc{block:StandardResult}{}
    with $ R_i = \langle \overline{z}_i \rangle $, and we get
    $$ Y \subseteq \prod_i \Sym ( \supp ( \overline{z}_i ) ) .$$
    By the hypothesis, $ Y $ is transitive, and therefore $ I $
    has size $ 1 $. Also $ z = \overline{z}_i $.
    This has the following consequences.
    \vspace*{1mm}
    \newline (i)   The orbit size must divide the order of
    $ \langle y \rangle $, that is $ s $ divides $ q $.
    \newline (ii)  The group $ \langle y, z \rangle $ is
    transitive on the support of $ Y $.
    \newline Since $ z $ and $ y $ commute, conjugation by $ z $
    also permutes the cyles of $ y $.
    \newline (iii) By (ii) $ z $ must be transitive on the
    cycles of $ y $ and then we have that the
    number of cycles $ m $ of $ y $ divides the order of
    $ \langle z \rangle $ which is $ p $. But $ p $ is prime,
    so either $ m = 1 $ or $ m = p $.    
    \newline {\large C}{\small ASE} 1    $ m = p $. Then we have $ {\rm d}1 ( Y ) = p q $.
    \newline {\large C}{\small ASE} 2    $ m = 1 $. Then $ {\rm d}1 ( Y ) = q = sp $.
    In this case, $ Y $ is a set of $ q $-cycles in $ G_Y \cong \Sym ( q ) $.    
}
\end {theorm}

\begin{block}
  { 
    We assume in this section that $ q $ is a prime.
    Theorem \arabicc{lemma:BoundOnIrreducibles}{} places a severe
    restriction on the possible degrees of the transitive fixed point
    sets.  We shall now compute the fixed point sets with these
    degrees. \par Firstly we consider the case $ p, q $ coprime. Thus
    Theorem \arabicc{lemma:BoundOnIrreducibles}{} asserts that any
    irreducible exact transitive fixed point set has degree $ q p $,
    which means we are
    dealing with the action of $ \Sym ( q p ) $ on its conjugacy class
    $ \Xi_{ q p } ^ q $ of fixed point free elements of order $ q $.
    Let $ X \subset \Xi_{ q p } ^ p $ be an
    irreducible exact transitive fixed point set
    and let $ x \in X $. Then a corresponding vertex $ Q_X $ fixes $ x
    $ and so $ Q_X \leq C_{ \Sym ( q p ) } ( x ) \cong {\mathbb Z}_q
    \wr \Sym ( p ) $. It follows that the size of $ Q_X $ is $ 1 $ or $
    p $, whence it must be $ p $ since $ X $ is transitive. Thus $ Q_X
    $ is generated by an element of order $ p $ which we denote by $ z
    $ in order to be consistent with the statement of Corollary
    \arabicc{corollary:StabP}{}. If $ x $ is written as the $ p $-fold
    product of $q$-cycles
    $$
    x = (a_{1 1} , a_{1 2}, \dots, a_{1 q })(a_{2 1} , a_{2 2},
    \dots, a_{2 q }) \dots (a_{p 1} , a_{p 2}, \dots, a_{p q }) $$
    then $ z $ may be taken to be the $ q $-fold product of $p
    $-cycles
    $$
    z = (a_{1 1} , a_{2 1}, \dots, a_{p 1 })(a_{1 2} , a_{2 2},
    \dots, a_{p 2 }) \dots (a_{1 q} , a_{2 q}, \dots, a_{p q }), $$
    and $ X = C_{ \Sym ( q p ) } ( Q_X ) \cap \Xi_{ q p } ^ q = C_{ {\rm Sym}
      ( q p ) } ( z ) \cap \Xi_{ q p } ^ q .$
        \par \vspace{5mm}
        Now we consider the case $ p = q $. In this case Theorem
        \arabicc{lemma:BoundOnIrreducibles}{} asserts that any
        transitive      fixed point set has degree $ p $ or $ p ^ 2 $.
        To find a transitive fixed point set $ X $ with $ d ( X ) = p
        $ we must consider      the action of $ \Sym ( p ) $ on its
        conjugacy class $ \Xi_{ p } ^ p $ of fixed point free elements
        of order $ p $. Now a   Sylow $ p $-subgroup of $ \Sym ( p ) $
        has order $ p $ so a corresponding vertex $ Q_X $ must have
        order $ p $, since      $ Q_X = 1 $ is ruled out by the
        transitivity of $ X $. Furthermore $ X = C_{ \Sym ( q p ) } (
        Q_X ) \cap \Xi_{ p } ^ p        = Q_X \cap \Xi_{ p } ^ p = Q_X
        \backslash 1 $.
        \par 
        To find a transitive fixed point set $ X $ with $ d ( X ) = p
        ^ 2 $ we must consider  the action of $ \Sym ( p ^ 2 ) $ on its
        conjugacy class $ \Xi_{ p^ 2 } ^ p $ of fixed point free
        elements of order $ p $.        Take $ x \in \Xi_{ p ^ 2 } ^ p
        $ and let $ C_x := C_{ \Sym ( p ^ 2 ) } ( x ) \cong {\mathbb
          Z}_p \wr \Sym ( p ) $. Now $ C_x $ contains    a Sylow $ p
        $-subgroup $ P $ of $ \Sym ( p ^ 2 ) $ and we may assume that $
        X $ is such that $ Q_X \leq P \leq C_x $, so $ x \in X $.
        By Corollary \arabicc{corollary:StabP}{} there is $ y \in X $
        such that $ Q_X = \Stab_P ( y ) $; hence $ Q_X $ fixes   both $
        x $ and $ y $, and also $ x \in Q_X $ since $ x \in \Stab_P ( y
        ) $. In the proof of Theorem
        \arabicc{lemma:BoundOnIrreducibles}{}   we used the fact that
        since $ y $ commutes with $ x $ it must induce a
        permutation of the $ p $ $ p $-cycles that compose $ x $, and
        by a remark in the third paragraph of that proof we see that
        this    induced action must be transitive, a condition on $ y
        $ that we wil refer to as condition (T).        It can be
        shown that for all $ y $ satisfying condition (T) we have
        $ \Stab_{ C_x } ( y ) = < x , y > $ and so $ Q_X \leq < x , y >
        $. Thus we either have  (i) $ Q_X = < x > $ or (ii) $ Q_X = <
        x , y > $. With (ii) it at first appears that different $ y $
        might yield non-conjugate       subgroups $ < x , y > $, but
        in fact it may be shown that any two elements $ y_1 $ and $
        y_2 $ in $ \Xi_{ p ^ 2 } ^ p $ satisfying condition (T) are
        conjugate in $ C_x $, and the resulting subgroups $ < x , y_1
        > $ and $ < x , y_2 > $ will be conjugate in $ C_x $. The
        corresponding   fixed point sets are easily computed.
        \par \vspace{5mm}
        We see that when $ p \neq q $ there is essentially at most one
        transitive fixed point set, and when $ p = q $ there are
        are most three essentially different fixed point sets. In
        either case the number of transitive fixed point sets is
        certainly finite.  }  { {\bf Transitive
        Fixed Point Sets For Prime $ q $} }
\end{block}
\begin{definition}
{
    Assume that $ H $ is any finite group and $ \Omega $ is any
    $ H $-set. Assume that the permutation module $ k \Omega $
    has a projective summand.
    \newline Define $ \kappa ( \Omega ) $ to be the lowest positive
    integer $ u $ such that the permutation module $ ( k \Omega ) ^ {\wr u} $
    for $ H \wr \Sym ( u ) $ does not have a projective summand, allowing
    the possibility $ \kappa ( \Omega ) = \infty $ if no such $ u $ exists.
    Recall \arabicc{lemma:WrProdHasProj}{} from which we get
    $  \kappa ( \Omega ) \geq p $. If $ j >  \kappa ( \Omega ) $ then
    $ k \Omega $ does not admit a projective summand by
    \arabicc{lemma:restrictedproductprojective}{}.
}
{    
}
\end{definition}
\begin{example}
{
    Let $ H = C_2 $, cyclic of order $ 2 $, and let $ p = 3 $.
    \newline (a) Take $ M = k \Omega $ with $ \Omega = H $ where
    $ H $ acts by left multiplication. Then consider $ M ^ {\wr 3} $,
    a module for $ H \wr \Sym ( 3 ) $.
    \newline If $ \langle \tau \rangle = H $ then take the basis
    $ \{ 1, \tau \} $ of $ M $, and take the basis for $ M ^ {\otimes 3} $
    consisting of all $ v_1 \otimes v_2 \otimes v_3 $ with $ v_i = 1 $
    or $ \tau $. Take a $ 3 $-cycle $ g $ of $ \Sym ( 3 ) $. This
    has two orbits of size $ 3 $ on this basis
    and two orbits of size $ 1 $.
    It follows that $ M ^ {\wr 3} $
    must have a projective summand as a module for
    $ \langle g \rangle $. But $ \langle g \rangle $ is a Sylow
    $ p $-subgroup of $ G := H \wr \Sym ( 3 ) $ and we deduce that
    $ M ^ {\wr 3} $ must have a projective summand as a module
    for $ G $.
    \newline So $ \kappa ( \Omega ) = 2 = p - 1 $.
    \vspace*{1mm}
    \newline (b) Now take $ M = k $, the trivial module, so that
    $ \Omega $ is a set with one element on which $ H $ acts
    trivially. Then $ M ^ {\wr a} $ is $ 1 $-dimensionsal
    for all $ a \geq 1 $ and hence does not have
    a projective summand for $ G $. So $ \kappa ( \Omega ) = \infty $.
}
{
}
\end{example}
\begin{definition}
{
     We shall call
     a set $ X \in  \mathcal S ^ q $ {\bf projective-free} if
     none of its irreducible components is projective. By
     \arabicc{block:Exact}{}
     products of exact sets are exact so projective-free
     implies exact.
}
{
}
\end{definition}
     \vspace*{1.5mm}
     \par
     Recall
     from \arabicc{corollary:FixIffProjComp}{} that if
     $ Z \in  \mathcal S ^ q $ is closed then $ Z $ is a
     fixed point set if and only if the $ k M_Z $-module
     $ k Z $ admits a projective summand.
\begin{theorm}
{
    Assume the setup of before, with elements of $ \mathcal S ^ q $.
    Assume $ X \in \mathcal S ^ q $.
    \newline (1) $ X $ is an irreducible exact fixed point set
    $ \Leftrightarrow $
    $ X = \Delta ^ { p^i } Y $ where $ Y $ is a transitive irreducible
    exact fixed point set and $ i \geq 0 $.
    \newline (2) $ X $ is a projective-free fixed point set
    $ \Leftrightarrow $
    $ X = X_1 ^ { a_1 } * \dots * X_t ^ { a_t } $ where the $ X_i $
    are pairwise coprime irreducible exact fixed point set,
    and $ 1 \leq a_i < \kappa ( X_i ) $.
    \newline (3) $ X $ is a fixed point set
    $ \Leftrightarrow $
    $ X = W * V $ where $ W $ is a projective-free fixed point set and $ V $
    is an irreducible projective fixed point set.
}
{
    (1) `$\Rightarrow$'
    By \arabicc{corollary:FormOfIrreducibles}{}
    $ X = \Delta ^ { p^i } Y $ for $ Y $ transitive
    (and therefore irreducible). From
    \arabicc{lemma::transitiverootclosed}{} $ Y $ is closed, so
    by \arabicc{projectivedeltanotclosed}{} it cannot be
    projective; and by \arabicc{lemma:ExactOrProjective}{}
    it must be exact.
    Also $ k Y $ admits a projective summand by
    \arabicc{lemma:StructureOfDeltaProduct}{}: so $ Y $ is
    a fixed point set.
    
    `$\Leftarrow$'
    By \arabicc{lemma:CoprimeSets}{} $ X $ is an irreducible
    fixed point set. The case $ i = 0 $ is trivial so assume
    $ i > 0 $. By \arabicc{lemma:StructureS}{} it is clear
    that $ Q_X \neq 1 $, so $ X $ is not projective and
    must be irreducible.

    (2) `$\Rightarrow$'
    Let $ X = X_1 ^ { a_1 } * \dots * X_t ^ { a_t } $
    be a decomposition into irreducibles.
    By \arabicc{theorm:ReductionToIrreducibles}{} each
    $ X_j $ is a fixed point set. But $ X_j ^ { a_j } $
    is also a fixed point set by 
    \arabicc{corollary:FactorsClosed}{}
    applied to $ X $  so
    by \arabicc{lemma:StructureOfProduct}{}
    $ 1 \leq a_j < \kappa ( X_j ) $. By assumption
    each $ X_j $ is not projective, hence exact.
    
    `$\Leftarrow$'
    We have $ X $ and each $  X_j ^ { a_j } $ is exact
    and closed by \arabicc{corollary:ProductIsClosed}{}.
    By assumption $ 1 \leq a_j < \kappa ( X_j ) $ so
    \arabicc{lemma:StructureOfProduct}{} tells us that
    each $ X_j ^ { a_j } $ is
    a fixed point set. It follows by
    \arabicc{lemma:CoprimeProductIsFixed}{} that
    $ X $ is a fixed point set.

    (3) `$\Rightarrow$'
    Factorize $ X $ into a product of irreducible sets. Say
    $$ X = W_1 ^ { a_ 1}  * \dots * W_s ^ { a_s } * V_1 ^ { b_ 1 } * \dots * V_s ^ { b_s } ,$$
    and write $ W = W_1 ^ { a_ 1}  * \dots * W_s ^ { a_s } $ and
    $ V = V_1 ^ { b_ 1 } * \dots * V_s ^ { b_s } $.
    Assume the $ W_i $ are exact and the $ V_j $ are projective.
    If there is no projective factor then there is nothing to do.
    So assume that $ t \geq 1 $.
    By
    \arabicc{corollary:FactorsClosed}{}
    $ V $ and $ W $ are fixed point sets.
    We wish to show that
    $ t = 1 $ and $ b = 1 $. 
    By
    \arabicc{corollary:FactorsClosed}{}
    each $ V_j ^ { b_j } $ is a fixed point set and 
    must be closed, which by
    \arabicc{lemma:projectiveproductnotclosed}{}
    means that $ b_j = 1 $. Now $ V $ is
    a fixed point set, hence closed, so by
    \arabicc{lemma:projectiveproductnotclosed}{}
    we get $ t = 1 $.            

    `$\Leftarrow$'
    We first show that $ X $ is closed.
    Since $ V $ is projective, we know
    that $ V = \Xi_V $ and because
    $ W $ and $ V $ are necessarily coprime
    we know by
    \arabicc{proposition:ProdBehaviour}{}
    that $ Q_X = Q_W * Q_V = Q_W * 1 $.
    The product $ W $ is projective-free, hence exact,
    which means
    that $ Q_W $ is fixed point free on $ \supp Q_W $.
    We can now compute directly that
    $$ {\rm Fix}_ { \Xi_ { W * V } } ( Q_W * 1 ) = W * V, $$
    and $ X $ is closed. 
    It follows from \arabicc{lemma:ClosedProductIsFixed}{}
    that $ X $ is a fixed point set.     
}
\end{theorm}

 \end{document}